\documentclass[11pt,reqno]{amsart}

\usepackage{latexsym}
\usepackage{amsfonts}
\usepackage{amsmath}
\usepackage{epsfig}
\usepackage{amssymb,amsmath}

\def\bl{\begin{lemma}}
\def\el{\end{lemma}}
\def\bth{\begin{theorem}}
\def\eth{\end{theorem}}
\def\bc{\begin{corollary}}
\def\ec{\end{corollary}}
\def\bcj{\begin{conjecture}}
\def\ecj{\end{conjecture}}
\def\bpr{\begin{proposition}}
\def\epr{\end{proposition}}
\def\bde{\begin{definition}}
\def\ede{\end{definition}}
\def\E{\mathbb{E}}

\newcommand{\be}{\begin{eqnarray}}
\newcommand{\ee}{\end{eqnarray}}
\newcommand{\eps}{{\mbox{$\epsilon$}}}

\newcommand{\A}{{\mathcal A}}
\newcommand{\B}{{\mathcal B}}

\newcommand{\RR}{{\mathcal R}}

\renewcommand{\and}{\hbox{ {\rm and} }}

\newcommand{\C}{{\mathcal{C}}}
\newcommand{\prob}{\mbox{\bf P}}
\newcommand{\p}{\mbox{\bf p}}

\newcommand{\lr}{\leftrightarrow}

\newcommand{\La}{{\mathcal{L}}}

\newtheorem{theorem}{Theorem}
\newtheorem{definition}{Definition}[section]
\newtheorem{lemma}[theorem]{Lemma}
\newtheorem{prop}[theorem]{Proposition}
\newtheorem{cor}[theorem]{Corollary}
\newtheorem{corollary}[theorem]{Corollary}
\newtheorem{proposition}[theorem]{Proposition}
\newtheorem{conjecture}[theorem]{Conjecture}

\theoremstyle{definition}
\numberwithin{equation}{section}

\input epsf.sty

\begin{document}
%\title{Necessary and sufficient conditions for mean-field behavior in ${1 \over
%d-1}$-bond-percolation}
\begin{center}
{{\Large{\bf{Mean-field conditions for percolation on finite
graphs}}}}
\end{center}

\author{{\sc Asaf Nachmias}}
\title{}
\thanks{Department of Mathematics, U.C. Berkeley and Microsoft Research. Research
supported in part by NSF grant \#DMS-0605166.}

\begin{abstract}
Let $\{G_n\}$ be a sequence of finite transitive graphs with
vertex degree $d=d(n)$
%(which may increase with $n$ or be bounded)
and $|G_n|=n$. Denote by $\p^t(v,v)$ the return probability after
$t$ steps of the {\em non-backtracking} random walk on $G_n$.
%In the spirit of Chung, Graham and Wilson \cite{CGW},
We show that if $\p^t(v,v)$ has {\em quasi-random} properties,
then critical bond-percolation on $G_n$ behaves as it would on a
random graph. More precisely, if
$$ \limsup_n n^{1/3} \sum _{t=1}^{n^{1/3}} t \p^t(v,v) < \infty \,
,$$ then the size of the largest component in $p$-bond-percolation
with $p={1 +O(n^{-1/3}) \over d-1}$ is roughly $n^{2/3}$. In
Physics jargon, this condition implies that there exists a scaling
window with a {\em mean-field} width of $n^{-1/3}$ around the
critical probability $p_c = {1 \over d-1}$.

%A slight variant of this condition also implies that above this
%scaling window the largest component is asymptotically larger than
%$n^{2/3}$.
A consequence of our theorems is that if $\{G_n\}$ is a transitive
expander family with girth at least $({2 \over 3} + \eps )
\log_{d-1} n$
%$g\geq {2 \over 3}
%\log_{d-1} n$
%$$ \limsup_n \Big \lfloor {g \over 2} \Big \rfloor - {1
%\over 3} \log_{d-1} n \geq 0 \, ,$$
then $\{G_n\}$ has the above scaling window around $p_c = {1 \over
d-1}$. In particular, bond-percolation on the celebrated Ramanujan
graph constructed by Lubotzky, Phillips and Sarnak \cite{LPS} has
the above scaling window. This provides the first examples of {\em
quasi-random} graphs behaving like random graphs with respect to
critical bond-percolation.

%bounded degree graph in which the scaling window assumes its
%mean-field width of $\Theta(n^{-1/3})$.

%This corollary also applies for percolation on transitive graphs
%with unbounded degree. In particular, we show that two-fold and
%three-fold products of complete graphs have a mean-field scaling
%window. This extends the work of van der Hofstad and Luczak
%\cite{HL}.
%In the realm of transitive expander graph with unbounded degree,
%our results show that any expander with  our results also apply
%for percolation on two-fold and three-fold products of complete
%graphs, extending the work of van der Hofstad and Luczak
\end{abstract}
\maketitle
\section{{\bf  \large Introduction}}

\subsection{ Background }Let $G$ be a graph and $p \in [0,1]$. Write $G_p$ for the graph
obtained from $G$ by performing $p$-bond-percolation on $G$, that
is, delete each edge with probability $1-p$ and retain it with
probability $p$, independently for all edges. Denote by $\C_1$ the
largest connected component of $G_p$.
%We call deleted edges {\em closed} and retained edges {\em open}.
When $G$ is the complete graph $K_n$, this model is known as the
Erd\H{o}s-R\'enyi random graph $G(n,p)$. Erd\H{o}s and R\'enyi
\cite{ER} discovered at 1960 that when $p_c = {1 \over n}$ the
model exhibits a phase transition. Namely, if $p = {c \over n}$
with $c<1$, then $|\C_1|$ is of order $\log n$ and if $c>1$, then
$|\C_1|$ is of order $n$ and all other components are of
logarithmic size.

The study of the random graph around the critical probability
(i.e., when $p \sim {1 \over n}$) was initiated by Bollob\' as
\cite{B1} over twenty years later. He showed that if $p = {1 +
O(n^{-1/3}) \over n}$ then the size of the largest component in
$G_p$ is roughly $n^{2/3}$. He also proved that if $p = {1 +
\eps(n) \over n}$ where $n^{1/3} \eps(n) \to \infty$ then with
high probability $n^{-2/3}|\C_1| \to \infty$, and if $n^{1/3}
\eps(n) \to -\infty$ then with high probability $n^{-2/3}|\C_1|
\to 0$ (Bollob\' as proved this with some logarithmic corrections
which were removed later by \L uczak \cite{Lu}). In Physics
jargon, this phenomenon is frequently called a {\em scaling
window} with {\em mean-field} width of $n^{-1/3}$ around $p_c = {1
\over n}$.

A scaling window of width $n^{-1/3}$ around $p_c={1 \over d-1}$
occurs also when $G$ is a random $d$-regular graph on $n$ vertices
where $d$ is fixed and $n \to \infty$ (see \cite{NP2} and
\cite{P2}). It is natural to expect, in the spirit of Chung,
Graham and Wilson \cite{CGW}, that critical percolation on {\em
deterministic} graphs having {\em quasi-random} properties will
behave the same as random graphs. Up to now, however, no examples
of this were known (see section \ref{related} for some related
results).

%of deterministic bounded degree graphs $G$ on which bond
%percolation has such a scaling window (in \cite{BCHSS1} it is
%proved that if $G$ is a high-dimension discrete torus, then a
%scaling window of size {\em at least} $n^{-1/3}$ exists, with no
%matching upper bound; see further discussion below).

\subsection{Mean-field scaling window} In this paper we show that if the
{\em non-backtracking} random walk (a simple random walk
restricted not to traverse the edge it has just visited in the
reverse direction. See section \ref{rw} for a precise definition)
behaves on $G$ as it would on a random graph in some sense, then
bond-percolation on $G$ has the same scaling window as it would on
the complete graph (or as it would on a random $d$-regular graph).
We give a quasi-random condition which guarantees the existence of
such a scaling window around $p_c = {1 \over d-1}$. This condition
is defined in terms of return probabilities $\p^t(v,v)$ of the
non-backtracking random walk and should be regarded as a {\em
geometric} condition on $G$.

\begin{theorem}\label{mainthm1} Let $\{ G_n\}$ be a family of
transitive graphs with vertex degree $d(n)\geq 3$ and assume
for simplicity that $|G_n|=n$. Let $p={1 + \lambda
n^{-1/3} \over d(n)-1}$ for some fixed $\lambda \in \mathbb{R}$
and consider the largest component $\C_1=\C_1(n)$ of
$p$-bond-percolation on $G_n$. Denote by $\p^t(v,v)$ the
probability that a non-backtracking random walk on $G_n$ starting
at a vertex $v$ will visit $v$ at time $t$. If \be
\label{maincond1} \limsup _{n} \, n^{1/3} \sum _{t=1}^{n^{1/3}} t
\p^t (v,v) < \infty \, , \ee then for any $\eps>0$ there exists
$A=A(\eps, \lambda)<\infty$ such that for all $n$
$$ \prob \Big ( {n^{2/3} \over A} \leq |\C_1| \leq An^{2/3}
\Big )  \geq 1 - \eps \, .$$
\end{theorem}

Condition (\ref{maincond1}) holds in various examples. In
particular, we have the following consequence for transitive
expander graphs. Recall that a sequence of connected graphs
$\{G_n\}$ is called an {\em expander} family if the largest
eigenvalue in absolute value, which is not $\pm 1$, of the
transition matrix of the simple random walk on $G_n$ is strictly
smaller than $1$, uniformly in $n$.

\begin{theorem}\label{expanderthm} Let $\{G_n\}$ be a transitive expander family with vertex degree $d(n)\geq 3$ and
assume that the girth $g(n)$ (i.e., the length of the shortest
cycle) of $G_n$ satisfies \be \label{expcond1} \limsup _n \,\,\Big
({ 1 \over d(n) -1} \Big ) ^{ \lfloor {g(n) \over 2} \rfloor}
n^{1/3} \log^2 n < \infty \, . \ee Then condition
(\ref{maincond1}) holds and in particular, if $\C_1$ is the
largest component of $p$-bond-percolation on $G_n$ with $p = {1 +
\lambda n^{-1/3} \over d(n)-1}$, then for any $\eps>0$ there
exists $A=A(\eps, \lambda)<\infty$ such that
$$ \prob \Big ( {n^{2/3} \over A} \leq |\C_1| \leq An^{2/3}
\Big )  \geq 1 - \eps \, .$$
\end{theorem}

%and $\eps(n)>0$ be a sequence satisfying $\eps(n)=o(1)$ but
%$\eps(n) n^{1/3} \to \infty$. We have
%\begin{enumerate}
%\item[(a)] { \rm (\bf Subcritical regime)} If $p={1 - \eps(n)
%\over d-1}$, then for any small $\delta >0$ we have
%$$ \prob \Big ( |\C_1| \geq \delta n^{2/3} \Big ) \longrightarrow
%0 \, , \qquad \hbox{\rm{as} } n \to \infty \, .$$
%
%\item[(b)] { \rm (\bf Critical window)} If
%
%\item[(c)] { \rm (\bf Supercritical regime)} If $p={1 + \eps(n)
%\over d-1}$, then for any large $A>0$ we have
%$$ \prob \Big ( |\C_1| \leq
%An^{2/3} \Big ) \longrightarrow 0  \, , \qquad\hbox{\rm{as} } n
%\to \infty \, . $$
%\end{enumerate}
%\end{theorem}

In particular, bond-percolation on the Ramanujan graphs
constructed by Lubotzky, Phillips and Sarnak \cite{LPS} have a
scaling window with mean-field width of $n^{-1/3}$ (the theorem
also applies for constructions of Margulis \cite{Mar}). These are
bounded degree expander graphs with girth approximately ${4 \over
3} \log_{d-1} n$, clearly satisfying the assumption of the theorem
with room to spare. Thus, Theorem \ref{expanderthm} gives the
first class of examples of {\em quasi-random} graphs having the
same behavior as random graphs with respect to critical
bond-percolation. We remark that it is proved in
\cite{KN} that the condition on the girth in Theorem \ref{expanderthm} is sharp up to the $\log^2 n$ factor. \\

%Theorem \ref{expanderthm} gives a class of examples of
%pseudo-random graphs which
%It is worth noting that part (a) of the theorem, and the upper
%bound on $|\C_1|$ of part (b) holds for {\em any} $G$ with maximum
%degree $d$, (see Proposition $1$ of \cite{NP2}). Thus, the novelty
%of Theorem~\ref{expanderthm} are the lower bounds on $|\C_1|$ of
%part (b) and (c).

Lower bounds on the size of percolation clusters at criticality is
a difficult task in general and is believed to depend heavily on
the geometry of the underlying graph. In the case $p_c={1 \over
d-1}$ overcoming this difficulty, not only establishes the
existence of a mean-field scaling window, but one immediately
reaps additional rewards: other geometric quantities of the
largest component, namely the diameter and the mixing time of the
simple random walk, assume their mean-field universal values.
Indeed, it is proved in \cite{NP3} that when $G$ has maximum
degree $d\in[3,n-1]$ and $p \leq {1 + O(n^{-1/3})\over d-1}$ {\em
if} $G_p$ typically has components of order $n^{2/3}$, then with
high probability, these components have diameter of order
$n^{1/3}$ and mixing time of order $n$. Hence, the following
corollary is an immediate consequence of Theorem \ref{mainthm1}
and of Theorem $1.2$ of \cite{NP3}.

\begin{cor}\label{mixdiam} Assume the setting of Theorem \ref{mainthm1}
and that condition (\ref{maincond1}) holds. Denote by {\rm
diam}$(\C_1)$ and by $T_{{\rm mix}}(\C_1)$ the diameter (maximal
graph distance) of $\C_1$ and the mixing time of the lazy simple
random walk on $\C_1$, respectively (see \cite{NP3} for a
definition). If $p={1 + \lambda n^{-1/3} \over d(n)-1}$ then for
any $\eps>0$ there exists $A=A(\eps, \lambda)$ such that
\begin{itemize}
\item $\displaystyle \qquad \prob \Big ( \hbox{\rm diam}(\C_1)
\not \in [A^{-1} n^{1/3}, An^{1/3}] \Big ) < \eps \, ,$

\item $\displaystyle \qquad \prob \Big ( T_{{\rm mix}}(\C_1) \not
\in [A^{-1} n, An] \Big ) < \eps \, .$ \\
\end{itemize}
\end{cor}
%\noindent In fact, one can prove that these large components of
%${1 \over d-1}$-bond-percolation are, in some sense, close to
%random trees on $n^{2/3}$ vertices (see \cite{NP3}). \\

\subsection{Outside the scaling window} In order to show that the
critical scaling window is of width $\Theta(n^{-1/3})$ around $p_c
= {1 \over d-1}$ one must show that if $p={1 + \eps(n) \over d(n)-1}$ and $|n^{1/3} \eps(n)| \to
\infty$, then $n^{-2/3} |\C_1|\to \infty$ with probability tending to $1$. The lower side of the window is easier to handle.
Indeed, in {\em any} $d$-regular graph, if $p={1 - \eps(n) \over d(n)-1}$ with $|n^{1/3} \eps(n)| \to
\infty$ (i.e., the subcritical regime), then $n^{-2/3} |\C_1| \to 0$ with probability tending to $1$. This is the contents
of part $(1)$ of Proposition $1$ of \cite{NP2}. Thus, we only need to take care of
the supercritical regime, $\eps(n) > 0$.

The next theorem shows that a slight variant of condition
(\ref{maincond1}) guarantees that $n^{-2/3} |\C_1| \to \infty$
with high probability when $p$ is above the mean-field scaling
window.

\begin{theorem} \label{mainthm2} Let $\eps(n)>0$ be a sequence such that
$\eps(n)=o(1)$ but $\eps(n) n^{1/3} \to \infty$. Take $p={1 +
\eps(n) \over d(n)-1}$ and $r=\eps^{-1} [ \log(n
\eps^3)-3 \log\log (n\eps^3)]$. Assume the setting of Theorem
\ref{mainthm1}. We have that if \be \label{maincond2} \limsup _{n}
\, \eps^{-1} r \sum _{t=1}^{2r} [(1+\eps)^{t\wedge r} -1] \p^t
(v,v) =0  \, , \ee then there is some fixed $\delta>0$ such that
$$ \prob \Big ( |\C_1| \geq {\delta \eps n \over \log^3 (n\eps^3) }\Big ) \longrightarrow 1 \,
, \quad \hbox{\rm as $n \to \infty$} \, ,$$ (note that
$ \delta \eps n \log^{-3}(n\eps^3) \gg n^{2/3}$).
\end{theorem}
\vspace{.2cm}

Condition (\ref{maincond2}) hold in various examples. Again, we
address expander graphs (even though conditions (\ref{maincond1})
and (\ref{maincond2}) do not require the graphs to be connected).

\begin{theorem} \label{expanderthm2} Under the setting of Theorem
\ref{expanderthm} we have that condition (\ref{maincond2}) holds.
In particular, if $p={1 + \eps(n) \over d(n)-1}$ with
$\eps(n)=o(1)$ but $\eps(n) n^{1/3} \to \infty$, then there is
some fixed $\delta>0$ such that
$$ \prob \Big ( |\C_1| \geq {\delta \eps n \over \log (n\eps^3) } \Big ) \longrightarrow 1 \,
, \quad \hbox{\rm as $n \to \infty$} \, .$$
\end{theorem}
\vspace{.2cm}

\subsection{Products of complete graphs} In \cite{HL}, van
der Hofstad and Luczak investigate bond-percolation on the
$k$-dimensional Hamming graph $H(k,m)$ which has vertex set
$\{0,\ldots, m-1 \}^k$ and a pair of vertices are connected by an
edge if and only if these vertices differ in precisely one
coordinate (this is a weak product of $k$ complete graphs, each on
$m$ vertices). In these graphs we have that $n= m^k$ and
$d(n)=k(m-1)$. It is shown in \cite{BCHSS1} and \cite{BCHSS2}, as
noted in \cite{HL}, that bond-percolation on $H(2,m)$ and $H(3,m)$
has a scaling window around $p_c={1 \over d-1}$ of width {\em at
least} $n^{-1/3}$ (i.e., the consequence of Theorem \ref{mainthm1}
hold). In \cite{HL} the authors provide an upper bound of order
$\log^{1/3}(n) n^{-1/3}$ on the width of the scaling window, i.e.,
they show that when $p={1 + \eps(n) \over d(n)-1}$ with $\eps(n) =
\Omega(\log^{1/3}(n) n^{-1/3})$ then $n^{-2/3}|\C_1| \to \infty$
with high probability. In fact, they show that for such $p$, the
size of the largest component is concentrated around $2 \eps(n)
n$.

In the following Theorem we prove that for $H(2,m)$ and $H(3,m)$
the scaling window is of width $\Theta(n^{-1/3})$ around $p_c = {1
\over d(n)-1}$. We remark that the same conclusion cannot be drawn
for $H(k,m)$ with $k>3$ by the results of \cite{KN}.
\begin{theorem}\label{hamming}
Conditions (\ref{maincond1}) and (\ref{maincond2}) hold for the sequences $H(2,m)$ and $H(3,m)$.
\end{theorem}

\subsection{Related results}\label{related}

It is worth comparing this work to results of a similar flavor
given by Borgs, Chayes, van der Hofstad, Slade and Spencer in
\cite{BCHSS1} and \cite{BCHSS2}. They define a {\em finite}
version of the Aizenman-Newman \cite{AN} triangle condition on
$\{G_n\}$ and show that it implies the existence of a scaling
window of width {\em at least} $n^{-1/3}$ around some $p_c$, in
which the largest component is of size about $n^{2/3}$. They
continue in \cite{BCHSS2} to show that this finite triangle
condition holds in a class of graphs including the
high-dimensional discrete finite torus $[m]^d$ (for large fixed
$d$ and $m \to \infty$) and the hypercube $\{0,1\}^{m}$.

The important advantage of the results of \cite{BCHSS1} over
Theorem \ref{mainthm1} is that it gives a mean-field scaling
window around a $p_c$ which is {\em not} necessarily ${1 \over
d-1}$. Indeed, in both the tori and the hypercube, at $p = {1 +
O(n^{-1/3}) \over d-1}$ we have that $n^{-2/3}|\C_1|\to 0$ with
high probability. Thus, the scaling window given in \cite{BCHSS1}
is at a higher location than ${1 + O(n^{-1/3}) \over d-1}$.
However, the advantage of Theorem \ref{maincond1} over the results
of \cite{BCHSS1} is that condition (\ref{maincond1}), if it holds,
is usually easy to verify while verifying the triangle condition
is notoriously difficult (even in the infinite case, see
\cite{HS}).

Another difference between the two results is that there is no
analogue in \cite{BCHSS1} to Theorem \ref{mainthm2}. Namely, it is
not known whether the finite triangle condition implies that
$n^{-2/3}|\C_1|\to \infty$ with high probability if $p=p_c(1 +
\eps(n))$ with $\eps(n)$ as in Theorem \ref{mainthm2}. Hence, the
finite triangle condition implies only that the scaling window is
of size $\Omega(n^{-1/3})$ but not $\Theta(n^{-1/3})$. In our
case, conditions (\ref{maincond1}) and (\ref{maincond2}) show that
the scaling window assumes the mean-field width of
$\Theta(n^{-1/3})$.

\subsection{Proof idea} Our proof is based on the analysis of the BFS (breadth-first-search)
exploration process on the percolated graph $G_p$ starting at a
uniformly chosen vertex $v$. It is obvious that this process is
``dominated'' by a BFS process on a percolated $d$-regular
infinite tree $T_p$. For instance, the size of the component
containing $v$ in $G_p$ is stochastically dominated by the size of
the component containing the root of $T_p$. This observation can be used to obtain
that the upper bound on $|\C_1|$ of Theorem \ref{mainthm1} (i.e.,
$|\C_1|$ is no more than $n^{2/3}$) holds for {\em any}
$d$-regular graph (this is shown in Theorem $1.2$ of \cite{NP3} or
Proposition $1$ of \cite{NP2}).

\begin{figure}\label{meanfieldpic}
\centering
\includegraphics[height = 6cm, width=15cm] {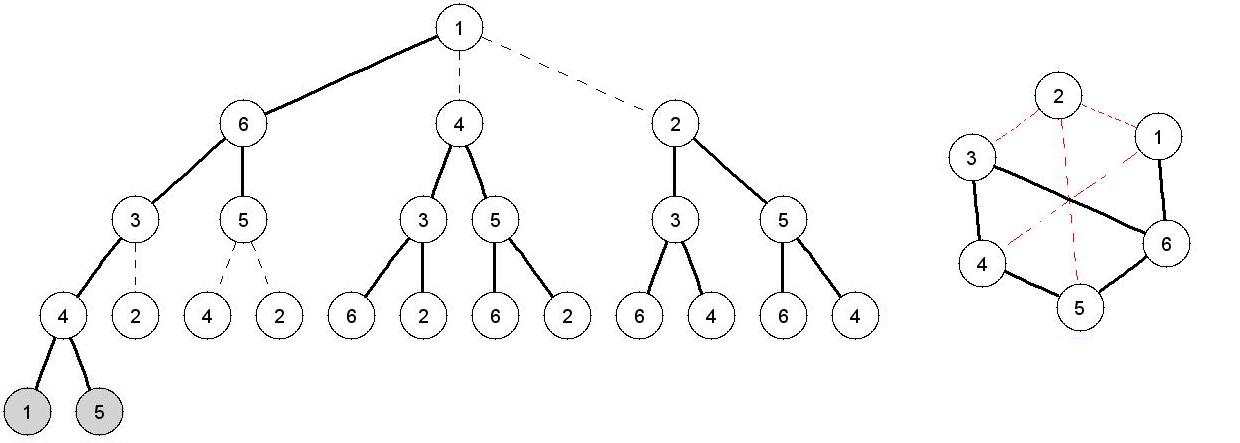}
\caption{On the left, the first levels of $T_p$ and on the right
$G_p$. Gray vertices on the left, labeled $5$ and $1$, are impure.}
\end{figure}

Obtaining a lower bound is much harder and requires additional
assumptions on the geometry of the underlying graph (e.g., the
triangle condition or condition (\ref{maincond1})). The starting
point of the approach of this paper is to consider how this
natural coupling of $G_p$ with $T_p$ fails to be sharp. This happens when we
encounter a vertex in $G_p$ which explores {\em less} than $d-1$
of its neighbors. This happens because at least one of its
neighbors has been explored before in the BFS process. We
consider $T$ as the covering tree of $G$ (i.e., to each vertex of
$T$ we associate a vertex of $G$ such that neighborhoods of
vertices of $G$ are preserved, see section $4.1$) and it is clear
that the above coupling fails when we explore, this time in $T_p$,
a vertex for which a vertex with the same label has been
discovered before, see figure $1$ (its full meaning will be clear in section $4.1$). We call such vertices {\em
impure} and a lower bound on the component of $G_p$ is given by
the component of $T_p$ after removing all impure vertices (and
their descendants). This is the contents of Proposition $11$.

The next important observation is that in some rough sense, the
unique path in $T$ between an impure vertex and the vertex causing
it to be impure (i.e., the vertex with the same label discovered
previously) is a random path in $T$ due to the nature of the BFS
exploration process. This path translates to a non-backtracking
random path in $G$. This enables us to bound from above the number
of impure vertices using our knowledge of the behavior of the
non-backtracking random walk on $G$. \\

This technique of estimating component sizes is novel. Previously
such estimates were obtained using the {\em lace expansion} (see
\cite{BCHSS1, BCHSS2}) and {\em sprinkling} (see \cite{BCHSS3,
HL}; sprinkling was introduced in \cite{AKS}).

%Hence, it is natural to try an achieve a similar lower bound on the component size. However, one cannot %expect this upper bound to be tight for all graphs. Indeed, it is easy to see that transitive graphs with %short cycles have $|\C_1|=O(\log n)$ when $p={1 \over d-1}$ and thus some additional geometry of the graph %is required for such a $n^{2/3}$ lower bound on $|\C_1|$ to hold.
%Lower bounds on component sizes in critical percolation tend to be significantly harder to obtain than %upper bounds. The two methods used before to this aim are the {\em lace expansion} in

%The novel approach in the proofs of this paper is obtaining a lower bound on component sizes by coupling %with percolation on a regular tree. The argument works roughly as follows. Imagine performing a BFS %process level by level on $G_p$ and coupling it naturally with percolation on a $d$-regular tree $T$, as %before. In the tree, each explored vertex explores $d-1$ new vertices while in $G$ it can happen that an %explored vertex has less than $d-1$ new vertices to explore because one of its neighbors have already been %seen before. The vertex of the tree corresponding to this neighbor has probability ${1 \over d-1}$ of %being connected to its father vertex and when this happens we obtain a vertex in $T_p$ which does not %exist in $G_p$. The crucial observation of this paper is that the path between this ``bad'' vertex of the %tree and vertex which makes it ``bad'' (i.e., the vertex of the tree corresponding

\subsection{Organization}
The rest of the paper is organized as follows. Deriving Theorems
\ref{expanderthm}, \ref{expanderthm2} and \ref{hamming} from
Theorems \ref{mainthm1} and \ref{mainthm2} is easy and done in the
next section. We present some useful preliminaries about
non-backtracking random walks and classical bond-percolation on
trees in Section \ref{pre}. Since the statements in Section \ref{pre}
are easy and classical we advise the reader to treat them as
black boxes and continue directly to Section \ref{lower}, which contains
the novel ideas of the proof. This section provides a coupling and a key lemma
allowing us to bound from below $|\C_1|$. We use the lemmas in
Section \ref{lower} in a straightforward manner to prove Theorems
\ref{mainthm1} and \ref{mainthm2}. We end with some concluding
remarks and open problems in Section \ref{conc}.

%Theorem \ref{mainthm1} and \ref{mainthm2} do not give any
%information when the critical percolation probability is above ${1
%\over d-1}$.

% Let $\{G_n\}$ be a transitive family of graphs, and
%for a vertex $v$ denote by $\C(v)$ the connected component
%containing $v$ in $p$-bond-percolation on $G_n$. The {\em
%susceptibility} $\chi(p)$ of $\{G_n\}$ is defined by
%$$ \chi (p) = \E |\C(v)| \, .$$ The
%authors of \cite{BCHSS1} {\em define} the critical probability
%$p_c$ of $\{G_n\}$ as the unique solution to the equation \be
%\label{pcdef} \chi(p_c) = \lambda n^{1/3} \, ,\ee for some small
%fixed $\lambda>0$ (this definition is motivated by the fact the
%susceptibility at the critical probability $p_c = {1 \over n}$ of
%the complete graph is indeed of order $n^{1/3}$). As in Theorem
%\ref{mainthm1}, the authors of \cite{BCHSS1} define a geometric
%condition on $\{G_n\}$, in particular, it is a {\em finite}
%version of the Aizenman-Newman triangle condition. In
%\cite{BCHSS1} they prove an analogue to Theorem \ref{mainthm1},
%namely, that the triangle condition implies the existence of a
%scaling window of size at least $n^{-1/3}$ around the $p_c$
%defined in (\ref{pcdef}).

%Verifying that the triangle condition holds is notoriously
%difficult. Using the lace expansion, the authors of \cite{BCHSS1}
%continued in \cite{BCHSS2} to prove that the triangle condition
%hold for the high-dimensional finite discrete torus $[n]^d$, the
%$n$-cube and some other examples. We note that there is no
%analogue of Theorem \ref{mainthm2} in \cite{BCHSS1}. Indeed, for
%the case

\section{ {\bf \large Percolation on expanders}} \label{percexp}

In this short section we derive Theorems \ref{expanderthm},
\ref{expanderthm2} and \ref{hamming} from Theorems \ref{mainthm1}
and \ref{mainthm2}. We refer the reader to Section \ref{rw} for a
formal definition of the non-backtracking random walk.

% To that aim we need to bound from above $\p^t(v,v)$ on expander
% graphs with girth $g$.  It is easy to see that if $\{G_n\}$ is an
% expander family, then the non-backtracking random walk Markov
% chain has a uniformly positive Cheeger constant. The results of
% Fill \cite{F} and Mihail \cite{Mi} show that even though our
% Markov chain is not reversible, rapid mixing does occur. It is an
% immediate deduction from these results that there exist some
% constant $c>0$ such for any set of vertices $A \subset G_n$ and
% for any vertex $v$ we have \be\label{mixing} \Big | \p^t(v, A) -
% {|A| \over n} \Big | \leq e^{-ct} \, ,\ee where $\p^t(v, A)$ is
% the probability that the
% non-backtracking random walk visits $A$ at time $t$. \\

\noindent {\bf Proof of Theorems \ref{expanderthm} and
\ref{expanderthm2}.}
%The results of Fill \cite{F} and Mihail
%\cite{Mi} show that even though the non-backtracking random walk
%is not a reversible Markov chain, rapid mixing does occur.
The proof simply requires to verify conditions (\ref{maincond1})
and (\ref{maincond2}). A recent result of Alon, Benjamini,
Lubetzky and Sodin \cite{ABLS} shows that the non-backtracking
random walk on $d$-regular, (with $d\geq 3$) non-bipartite
expanders mixes faster than the usual simple random walk. It
follows from their results that in that case there exists some
$C>0$ such that for all $t \geq C \log n$ we have $\p ^t(v,v) \leq
{2 \over n}$. If the expander graph $G$ happens to be bipartite
(as in the case of the Lubotzky-Phillips-Sarnak graph \cite{LPS}),
it is clear that $p^t(v,v)=0$ if $t$ is odd. Let $A$ and $B$ be
the two parts of the graph. We may consider the connected
component of $G^2$ induced on the vertices of $A$. One can readily
show that this graph is an expander on $n/2$ vertices, and thus
the results of \cite{ABLS} apply. We learn from this discussion
that if $G_n$ is an expander family (either bipartite or
non-bipartite), then there exists some fixed $C>0$ such that for
all $t \geq C \log n$ we have $\p ^t(v,v) \leq {4 \over n}$.

To handle smaller $t$'s we use the girth assumption (in the same
way done in \cite{ABLS}).
%Indeed, let
%$$A = \Big \{ u \in G : d_G(u,v) = \Big \lfloor {g \over 2} \Big \rfloor \Big \} \, ,$$
%where $d_G(u,v)$ is the usual graph distance between $u$ and $v$.
Since the graph spanned on $\{ u : d_G(u,v) \leq \lfloor g/2
\rfloor \}$ is a tree, in order for the walk to return to $v$ at
time $t$, it must visit the set $\{ u : d_G(u,v) = \lfloor g/2
\rfloor \}$ at time $t- \lfloor g/2 \rfloor$ and then take
precisely $\lfloor g/2 \rfloor$ steps towards $v$ in the tree.
Hence for any $t\geq g$ we have
$$ \p^t(v,v) \leq  \Big ( {1 \over d-1} \Big )^{\lfloor g/2
\rfloor} \, ,$$ and for $t < g$ it is clear that $\p^t(v,v) = 0$.
To sum things up, we have
$$ \p^t(v,v) \leq \left \{
\begin{array}{ll} 0\, , & \,\, t < g \, , \\
\Big ( {1 \over d-1} \Big )^{\lfloor g/2 \rfloor} \, , & \,\, g \leq t < C \log n \, , \\
{4 \over n} \, ,& \, \, t \geq C \log n \, . \\
\end{array} \right . $$
We use this to verify that condition (\ref{maincond1}) holds,
\begin{eqnarray*} n^{1/3} \sum _{t=1}^{n^{1/3}} t\p^t(v,v) &\leq& n^{1/3} \sum
_{t=g}^{C \log n} t \Big ( {1 \over d-1} \Big )^{\lfloor g/2
\rfloor} + n^{1/3} \sum_{t=C \log n}^{n^{1/3}} {4 t \over n} \\
&\leq& C^2 n^{1/3} \log^2 n \Big ( {1 \over d-1} \Big )^{\lfloor
g/2 \rfloor} + 2 \, .\end{eqnarray*} Thus, assumption
(\ref{expcond1}) on the girth implies that condition
(\ref{maincond1}) holds. We now verify that condition
(\ref{maincond2}) holds,
\begin{eqnarray} \label{calcend} \eps^{-1} r \sum _{t=1}^{2r} [(1+\eps)^{t \wedge r} -1] \p^t
(v,v) &\leq& {\eps^{-1} r \over (d-1)^{\lfloor g/2 \rfloor}}
\sum_{t=g}^{C \log n} [(1+\eps)^{t\wedge r} -1] \nonumber \\ &+&
{4 \eps^{-1} r \over n} \sum _{t=C \log n}^{2r} (1+\eps)^{t\wedge
r} \, .
\end{eqnarray} We estimate the second term on the right hand size
with
$$ \sum _{t=C \log n}^{2r} (1+\eps)^{t\wedge r} \leq (\eps^{-1} + r)
(1+\eps)^r\, .$$  Recall that $r = \eps^{-1} [ \log(n\eps^3) - 3\log \log (n\eps^3)]$, hence the second term on the
right hand side of (\ref{calcend}) is of order
$$ \eps^{-1} n^{-1} r^2 (1+\eps)^r = O(\log ^{-1} (n\eps^3)) = o(1) \, ,$$ by our assumption on
$\eps$. To estimate the first term on the right hand size of
(\ref{calcend}), note that there exists $C_2>0$ such that for all
$t \leq \eps^{-1}$ we have $(1+\eps)^t \leq 1 + C_2 \eps t$. Thus,
in the case that $\eps^{-1} \geq C \log (n)$ we estimate this term
by
$$  {\eps^{-1} r \over (d-1)^{\lfloor g/2 \rfloor}}
\sum_{t=g}^{C \log n} [(1+\eps)^{t\wedge r} -1] \leq { O( r \log
^2 n ) \over (d-1)^{\lfloor g/2 \rfloor}} = o(1) \, ,$$ by
(\ref{expcond1}) (note that $r = o(n^{1/3})$). If $\eps^{-1} \leq
C \log (n)$ we estimate
$$ {\eps^{-1} r \over (d-1)^{\lfloor g/2 \rfloor}}
\sum_{t=g}^{C \log n} [(1+\eps)^{t\wedge r} -1] \leq {\eps^{-2} r
(1+\eps)^{C \log n} \over (d-1)^{\lfloor g/2 \rfloor}} \leq {
\eps^{-3} \log(n \eps^3) n^{C \eps} \over (d-1)^{\lfloor g/2
\rfloor}}  = o(1)\, ,$$ by our assumption on $\eps(n)$ and
(\ref{expcond1}). \qed \\

\noindent{\bf Proof of Theorem \ref{hamming}.} The graphs $H(2,m)$
and $H(3,m)$ are expanders with girth $3$ (i.e., they have
triangles). Assumption (\ref{expcond1}) can be seen to hold for
$H(2,m)$, however, it does not hold for $H(3,m)$ because of the
$\log^2 n$ term appearing in (\ref{expcond1}). Therefore, we need
to prove that conditions (\ref{maincond1}) and (\ref{maincond2})
hold for $H(3,m)$, and to that aim we estimate more carefully
$\p^t(v,v)$. Take $v=(0,0,0)$ and define the following subsets of
the vertex set of $H(3,m)$,
$$ A_1 = \Big \{ (i,j,k) \in H(3,m) : i=j=0 \hbox{ {\rm and} } k
\neq 0 \Big \} \, ,$$
$$ A_2 = \Big \{ (i,j,k) \in H(3,m) : i=k=0 \hbox{ {\rm and} } j
\neq 0 \Big \} \, ,$$
$$ A_3 = \Big \{ (i,j,k) \in H(3,m) : j=k=0 \hbox{ {\rm and} } i
\neq 0 \Big \} \, .$$ Observe that in order for the
non-backtracking random walk to return to $v$ at time $t$, it must
be in $A_1 \cup A_2 \cup A_3$ at time $t-1$. Hence,
$$ \p^t(v,v) \leq {1 \over d(n) - 1} \p^{t-1}(v, A_1 \cup A_2 \cup
A_3) \, ,$$ where $\p^{t-1}(v,A)$ is the probability that the walk
visits $A$ at time $t-1$. Let $\{X_t\}$ be the non-backtracking
random walk. We have \begin{eqnarray*} \p^{t-1}(v, A_1) &\leq&
\prob \Big (
X_j \in A_1 \cup \{v\} \hbox{ {\rm for all } } j \leq t-1 \Big )\\
&+& \prob \Big ( \exists j \leq t-1 \hbox{ {\rm with } } X_j \not
\in A_1 \cup \{v \} \hbox{ {\rm and} } X_{j+1} \in A_1 \Big ) \,
.\end{eqnarray*} It is clear that the probability of the first
event is $({2\over 3})^{t-1}$ because it requires that the walk
does not walk on coordinates $2$ and $3$ for $t-1$ steps. We can
bound the probability of the second event above by ${t \over
d(n)-1}$, since the probability of a particular $j$ having $X_j
\not \in A_1 \cup \{v \}$ and $X_{j+1} \in A_1$ is bounded above
by ${1 \over d(n)-1}$ (because it needs to walk to $0$ at the
first coordinate). We deduce by all this that
$$ \p^t(v,v) \leq \left \{
\begin{array}{ll} 0\, , & \,\, t < 3 \, , \\
 {3 \over d(n)-1} \Big ( ({2 \over 3})^{t-1} + {t \over d(n) -1} \Big ) \, , & \,\, g \leq t < C \log n \, , \\
{2 \over n} \, ,& \, \, t \geq C \log n \, . \\
\end{array} \right . $$
We use this to verify that condition (\ref{maincond1}) holds,
\begin{eqnarray*} n^{1/3} \sum _{t=1}^{n^{1/3}} t\p^t(v,v) &\leq& {3 n^{1/3} \over d(n)-1}  \sum
_{t=3}^{C \log n} t (2/3)^{t-1} + {3 n^{1/3} \over (d(n)-1)^2}
\sum _{t=3}^{C\log n} t^2 + n^{1/3} \sum_{t=C \log n}^{n^{1/3}} {2
t \over n}  \, .\end{eqnarray*} Recalling that $d(n) =
\Theta(n^{1/3})$ shows that condition (\ref{maincond1}) holds. We
now verify that condition (\ref{maincond2}) holds,

\begin{eqnarray*} \eps^{-1} r \sum _{t=1}^{2r} [(1+\eps)^{t \wedge r} -1] \p^t
(v,v) &\leq& {3\eps^{-1} r \over d(n)-1}
\sum_{t=3}^{C \log n} [(1+\eps)^{t\wedge r} -1] \Big ({2 \over 3} \Big )^{t-1} \nonumber \\
&+& {3\eps^{-1} r \over (d(n)-1)^2} \sum_{t=3}^{C\log
n}[(1+\eps)^{t\wedge r} -1] t \\ &+& {\eps^{-1} r \over n} \sum
_{t=C \log n}^{2r} (1+\eps)^{t\wedge r} \, .
\end{eqnarray*}
The last term on the right hand side tends to $0$ as in
(\ref{calcend}). To estimate the two other terms, assume first
$\eps^{-1}(n) \geq C \log n$, then as before $(1+ \eps)^t \leq 1 +
C_2 \eps t$ and we have
$$ {3\eps^{-1} r \over d(n)-1}
\sum_{t=3}^{C \log n} [(1+\eps)^{t\wedge r} -1] \Big ({2 \over 3}
\Big )^{t-1} \leq  {O(r) \over d(n) -1 } \sum_{t=3}^\infty t \Big(
{2 \over 3} \Big )^{t-1} = o(1) \, ,$$ since $r=o(n^{1/3})$.
Similarly,
$$ {3\eps^{-1} r \over (d(n)-1)^2} \sum_{t=3}^{C\log
n}[(1+\eps)^{t\wedge r} -1]t \leq { O(r \log^3 n) \over (d(n)-1)^2
} = o(1) \, .$$ The case $\eps^{-1}(n) \leq C \log n$ is handled
similarly and we conclude that conditions
(\ref{maincond1}) and (\ref{maincond2}) indeed hold. \qed \\

%We split the first sum on the right hand size of into to two sums
%in the domains $g \leq t \leq \eps^{-1}$ and $\eps^{-1} \leq t
%\leq 2r$. There exists a fixed $C>0$ such that for all $t \leq
%\eps^{-1}$ we have $(1+\eps)^t \leq 1 + C\eps t$, thus
%$$ {\eps^{-1} r \over (d-1)^{\lfloor g/2 \rfloor}}
%\sum_{t=g}^{\eps^{-1}} [(1+\eps)^t -1] e^{-c(t-\lfloor g/2
%\rfloor)} \leq {C r \over (d-1)^{\lfloor g/2
%\rfloor}}\sum_{t=g}^{\eps^{-1}} t e^{-c(t-\lfloor g/2 \rfloor)} =
%o(1) \, ,$$ by our assumption on the girth and our assumption on
%$\eps$ (note that $r = o(n^{1/3})$. We also have
%$$ {\eps^{-1} r \over (d-1)^{\lfloor g/2 \rfloor}}
%\sum_{t=\eps^{-1}}^{2r} [(1+\eps)^{t \wedge r} -1] e^{-c(t-\lfloor
%g/2 \rfloor)} \leq {\eps^{-1} r \over (d-1)^{\lfloor g/2 \rfloor}}
%\sum_{t=\eps^{-1}}^{\infty} e^{-ct/2} \leq { e^{-{c \eps^{-1}\over
%4}} \eps^{-1} r \over (d-1)^{\lfloor g/2 \rfloor}} \, ,
%$$
%which also tends to $0$ as $n\to \infty$ by our assumption on
%$\eps$, concluding our proof.
%\qed \\

%In \cite{ER} the authors studied bond percolation on the complete
%graph, a model known as the Erd\H{o}s-R\'enyi random-graph
%$G(n,p)$ which is obtained from the complete graph on $n$ vertices
%by retaining each edge with probability $p$ and deleting it with
%probability $1-p$, independently of all other edges. They proved
%that the model exhibits a phase transition when $p$ is scaled such
%that $p={c \over n}$ where $c>0$ is constant.

\section{{\bf \large Preliminaries}}\label{pre}

\subsection{The non-backtracking random walk}\label{rw} The non-backtracking random walk
is a simple random walk on a graph not allowed to traverse back on
an edge it has just walked on. Formally, the non-backtracking
random walk on an undirected graph $G=(V,E)$, starting from a
vertex $x \in V$, is a Markov chain $\{X_t\}$ with transition
matrix $\prob ^x$ on the state space of {\em directed} edges
$$ \overrightarrow{E} = \Big \{ (x,y) \, : \, \{x,y\} \in E \Big \} \, .$$
If $X_t = (x,y)$ we write $X_t^{(1)} = x$ and $X_t^{(2)} = y$.
Also, for notational convenience, we write $\prob _{(x,w)}(\cdot)$
for $\prob ^x( \cdot \mid X_0 = (x,w) )$, and $\p^t(x,y)$ for
$\prob^x (X_t^{(2)} = y)$. The non-backtracking walk starting from
a vertex $x$ has initial state given by
$$\prob ^x(X_0 = (x,y)) = {\bf 1}_{\{(x,y)\in \overrightarrow{E}\}} {1 \over {\rm deg}(x)} \, ,$$
and transition probabilities given by
$$ \prob^x_{(x,y)}(X_1 = (y,z)) = {\bf 1}_{\{(y,z)\in \overrightarrow{E} \, , z \neq
x \}} {1 \over {\rm deg}(y) - 1} \, ,$$ where we write deg$(x)$
for the degree of $x$ in $G$. The following Lemma is the analogue
of the statement $\p^{t+t'}(x,x) = \sum _{y} \p
^t(x,y)\p^{t'}(y,x)$ which holds for the simple random walk.

\begin{lemma}\label{loop} Let $G=(V,E)$ be a transitive graph
with vertex degree $d$ and $v\in V$ be an arbitrary vertex. Denote
$v$'s neighbors in $G$ by $v_1, \ldots, v_d$. Then for any two
positive integer $t, t'$ we have $$ \sum _{y \in V}
\mathop{\sum_{i,j=1}^{d}}_{i\neq j} \prob_{(v,v_i)}(X_t^{(2)}=y)
\prob_{(v,v_j)}(X_{t'}^{(2)}=y) = d(d-1) \p^{t+t'+1}(v,v) \, .$$
\end{lemma}
\noindent {\bf Proof.} We expand $\p^{t+t'+1}(v,v)$ by
conditioning on the location of the chain at time $t+1$. The
Markov property gives, \be \label{nonbtstep1} \p^{t+t'+1}(v,v) =
\sum _{ e = (x,y)} \prob_{(x,y)} \Big ( X_{t'}^{(2)} = v \Big )
\prob ^v \Big (X_{t+1} = (x,y)\Big ) \, .\ee For a vertex $x \in
V$ write $x_1, \ldots, x_d$ for the neighbors of $x$ in $G$. We
use
$$ \prob ^v \Big (X_{t+1} = (x,y)\Big ) = {1 \over d-1} \mathop{\sum_{j=1}^d}_{x_j \neq
y} \prob ^v \Big (X_t = (x_j,x)\Big ) \, ,$$ to rewrite
(\ref{nonbtstep1}) as \be \label{nonbtstep2} \p^{t+t'+1}(v,v) = {1
\over d-1} \sum _{x \in V} \mathop{\sum_{i,j=1}^{d}}_{i\neq j}
\prob_{(x,x_i)} \Big ( X_{t'}^{(2)} = v \Big ) \prob ^v \Big (X_t
= (x_j,x)\Big ) \, .\ee Write $N(v,x_j,x,t)$ for the number of
non-backtracking paths of length $t$ from $v$ to $x$ such that the
directed edge visited at time $t$ (the last edge) is $(x_j,x)$. We
have $\prob ^v \Big (X_t = (x_j,x)\Big ) =
N(v,x_j,x,t)/d(d-1)^{t-1}$. By traversing the paths in reverse we
learn that $N(v,x_j,x,t)$ is the number of non-backtracking paths
of length $t$, starting with the edge $(x,x_j)$ and ending in $v$.
Hence we have $\prob _{(x,x_j)} \Big ( X_t^{(2)} = v \Big ) =
N(v,x_j,x,t) / (d-1)^{t-1}$. We deduce that
$$ \prob ^v \Big (X_t = (x_j,x)\Big ) = {1 \over d} \prob _{(x,x_j)} \Big
( X_t^{(2)} = v \Big ) \, .$$ We put this into (\ref{nonbtstep2})
and get
$$ \p^{t+t'+1}(v,v) = {1
\over d(d-1)} \sum _{x \in V} \mathop{\sum_{i,j=1}^{d}}_{i\neq j}
\prob_{(x,x_i)} \Big ( X_{t'}^{(2)} = v \Big ) \prob _{(x,x_j)}
\Big ( X_t^{(2)} = v \Big ) \, .$$ We sum on $v\in V$ and divide
by $n$ both sides to get \be \label{nonbtstep3} {1 \over n} \sum
_{v \in V} \p^{t+t'+1}(v,v) = {1 \over d(d-1)n} \sum _{x \in V}
\sum _{v \in V} \mathop{\sum_{i,j=1}^{d}}_{i\neq j}
\prob_{(x,x_i)} \Big ( X_{t'}^{(2)} = v \Big ) \prob _{(x,x_j)}
\Big ( X_t^{(2)} = v \Big ) \, .\ee Observe that because $G$ is
transitive we have ${1 \over n} \sum _{v \in V} \p^{t+t'+1}(v,v) =
\p^{t+t'+1}(v,v)$. Also due to transitivity we have that for any
$x \in V$ the sum
$$ \sum _{v \in V}
\mathop{\sum_{i,j=1}^{d}}_{i\neq j} \prob_{(x,x_i)} \Big (
X_{t'}^{(2)} = v \Big ) \prob _{(x,x_j)} \Big ( X_t^{(2)} = v \Big
) \, ,$$ evaluates to the same number. The assertion of the lemma
then follows from (\ref{nonbtstep3}). \qed

\subsection{Critical percolation on trees}\label{trees} Let $T$ be an infinite
$d$-regular tree rooted at a vertex $\rho$. For a vertex $u \in T$
we write $|u|$ for the distance of $u$ from $\rho$ (i.e., the
number of edges in the path between $u$ and $\rho$). We say $w$ is
an {\em ancestor} of $u$ (or $u$ is a {\em descendant} of $w$) if
$w$ belongs to the unique path connecting $\rho$ and $u$. For a
vertex $w \in T$ with $w \neq \rho$, denote by $w^-$ the immediate
ancestor (the father) of $w$ in $T$. We also use the notation $u
\wedge w$ for the common ancestor of $u$ and $w$ having maximal
distance from $\rho$.

Let $p\in [0,1]$ and consider $p$-bond percolation on $T$. Denote
the resulting subgraph by $T_p$. For two vertices $u, w \in T$, we
denote by $\{ u \lr w \}$ the event that $u$ is connected to $w$
in $T_p$. For an integer $r>0$ we write $H_r$ for the $r$-th level
of $T_p$, i.e.,
$$ H_r = \Big | \Big \{ w \in T \, : \, |w|=r \, , w \lr \rho \Big \} \Big | \,
.$$ \\ \indent We will frequently use the following two lemmas
dealing with percolation on $T$ with $p={1 + \eps \over d-1}$. The
two lemmas estimate the same quantities, but since percolation at
this regime changes drastically with the sign of $\eps$, the
estimates given in them are different.

% Note
%that the estimates in this lemma for $\eps<0$ are quite crude (for
%$\eps>0$ they are of the correct magnitude), however, they will
%only effect the constants the proofs yield.

\begin{lemma} \label{treecalc} Let $\eps\in (0,1/2)$ and put $p={1 + \eps \over d-1}$.
For any integer $r>0$ we have \be \label{treemom2} \E H_r^2 = O
\Big ( \eps^{-1} (1+\eps)^{2r} \Big ) \, , \ee \be
\label{surviveprob} {\eps \over 2} \leq \prob \Big ( H_r
> 0 \Big ) \leq 12 \eps (1-e^{-\eps r/2})^{-1} \,
,\ee \be \label{treemom2cond} \E \Big [ (\sum _{k=r/2}^{r} H_k)^2
\mid H_{r/2}>0 \Big ] = O(\eps^{-4} (1+\eps)^{2r} \Big ) \, .\ee
\end{lemma}

\begin{lemma} \label{treecalcneg} Let $\eps\in (0,1/2)$ and put $p={1 - \eps \over d-1}$.
For any integer $r>0$ we have \be \label{treemom2neg} \E H_r^2 = O
\Big ( \eps^{-1} (1-\eps)^{r} \Big ) \, , \ee \be
\label{surviveprobneg} {\eps (1-\eps)^r \over 2} \leq \prob \Big (
H_r
> 0 \Big ) \leq {12 \over r} \,
,\ee

\be \label{treemom2condneg} \E \Big [ (\sum _{k=r/2}^{r} H_k)^2
\mid H_{r/2}>0 \Big ] = O( \eps^{-3} r ) \, .\ee
\end{lemma}

\vspace{.1cm} \noindent For the proof of (\ref{surviveprob}) and
(\ref{surviveprobneg}) we will use the following result due to
Lyons \cite{L}.

\begin{lemma} [Theorem 2.1 of \cite{L}] \label{lyons}
Assign each edge $e$ from level $b-1$  to level $b$ of $T$ the
edge resistance $r_e = {1-p \over p^b}$. Let $\RR_k$ be the
effective resistance from the root to level $k$ of $T$. Then
\begin{equation} \label{russ} { 1 \over 1 + \RR_k } \leq \prob ( H_k > 0 ) \, \leq { 2 \over 1 + \RR_k }  \, .
\nonumber \end{equation}
\end{lemma}
\vspace{0.1 cm}

\noindent {\bf Proof of Lemma \ref{treecalc}.} We have
$$ \E H_r^2 = \sum _{j=0}^{r} f(j) p^{2r-j} \, ,$$
where
$$f(j) = \Big |\Big \{ \{w_1, w_2\} : |w_1|=|w_2|=r \, , \, |w_1
\wedge w_2|=j  \Big \} \Big | \, .$$ We have
$f(0)=d(d-1)^{2r-1}/2$ and $f(r)=d(d-1)^{r-1}$ and $f(j) =
d(d-2)(d-1)^{2r-j-2}/2$ for $1 \leq j \leq r-1$. We get
\begin{eqnarray} \label{uselater} \E H_r^2  &=& {d (1+\eps)^{2r} \over
2(d-1)} + {d (1+\eps)^r \over d-1} + {d(d-2)(1+\eps)^{2r} \over
2(d-1)^2} \sum _{j=1}^{r-1} (1+\eps)^{-j}  \\ &=& O \Big (
{(1+\eps)^{2r} \over \eps} \Big ) \, , \nonumber \end{eqnarray}
%( (1+\eps)^{2k} {1 - (1- {\eps \over (1+\eps)})^{k} \over \eps}
%\Big ) = \Theta \Big ( { (1+\eps)^{2k} (1-e^{-\eps k}) \over \eps}
%\Big ) \, , \end{eqnarray*}
which finishes the proof of (\ref{treemom2}). To prove
(\ref{surviveprob}), in the setting of Lemma \ref{lyons}, we have
that the effective resistance $\RR_r$ from $\rho$ to level $r$ of
$T$ satisfies (see \cite{P}, Example 8.3)
\begin{eqnarray} \label{survstep} \RR_r  &=& \sum _{i=1}^r {(1-p) p^{-i} \over d(d-1)^{i-1}}
= \sum _{i=1}^r { {d-2-\eps \over d-1} (d-1)^{i} (1+\eps)^{-i}
\over d(d-1)^{i-1}} \nonumber \\ &=& {d-2-\eps \over d} \Big [
\eps^{-1} (1-(1+\eps)^{-r}) \Big ] \, .\end{eqnarray} We bound the
last term using the estimates
$$ {1 \over 6} \leq {d-2-\eps \over d} \leq 1  \, , \qquad
0 \leq (1+\eps)^{-r} \leq e^{-\eps r/2} \, ,$$ which are valid for
$d \geq 3$ and $\eps\in (0,1/2)$ (since $1+x \geq e^{x/2}$ for
$x\in[0,1/2]$). We get that
$$  {\eps^{-1} (1- e^{-\eps r/2})\over 6} \leq \RR_r \leq \eps^{-1}  \, ,$$
which together with Lemma \ref{lyons} yields (\ref{surviveprob}).
To prove (\ref{treemom2cond}) note that for any $k_1 \geq k_2 \geq
r/2$ we have
$$ \E \Big [ H_{k_1} H_{k_2} \mid H_{r/2}>0 \Big ] = {\E H_{k_1} H_{k_2}
\over \prob (H_{r/2}>0)} \, .$$ Thus, \be \label{mom2step1} \E
\Big [ (\sum _{k=r/2}^{r} H_k)^2 \mid H_{r/2}>0 \Big ] = { 2 \over
\prob (H_{r/2} > 0)}  \sum_{k_1 \geq k_2 \geq r/2}^r {\E H_{k_1}
H_{k_2} }\, .\ee For any $k_1 \geq k_2$ we have
$$ \E [ H_{k_1} H_{k_2} ] = \E [ H_{k_2} \E [H_{k_1} \mid H_{k_2}]
] = (1+\eps)^{k_1-k_2} \E[H_{k_2}^2] = O \Big ( \eps^{-1}
(1+\eps)^{k_1 + k_2} \Big ) \, ,$$ by (\ref{treemom2}). We put
this into (\ref{mom2step1}) and use the lower bound of
(\ref{surviveprob}) to get that $$\E \Big [ (\sum _{k=r/2}^{r}
H_k)^2 \mid H_{r/2}>0 \Big ] \leq O(\eps^{-2}) \sum _{k_1=r/2}^r
\sum _{k_2=r/2}^{k_1} (1+\eps)^{k_1 + k_2} = O(\eps^{-4}
(1+\eps)^{2r}
\Big ) \, ,$$ finishing the proof of (\ref{treemom2cond}). \qed \\

\noindent {\bf Proof of Lemma \ref{treecalcneg}.} We proceed in
the same manner as in the previous lemma. Indeed, equality
(\ref{uselater}) and an immediate calculation imply
(\ref{treemom2neg}). The equality (\ref{survstep}) on $\RR_r$
holds and we get
$$ \RR_r = {d-2+\eps \over d} \Big [ \eps^{-1} ( (1-\eps)^{-r}
-1) \Big ] \, .$$ We use $(1-\eps)^{-r}-1\geq \eps r$ and similar
estimates to previous lemma to get
$$ {r \over 6} \leq \RR_r \leq \eps^{-1}
(1-\eps)^{-r} \, ,$$ which together with Lemma \ref{lyons} yields
(\ref{surviveprobneg}). As before, using (\ref{treemom2neg}) we
get that for any $k_1 \geq k_2$
$$ \E H_{k_1} H_{k_2} = (1-\eps)^{k_1 - k_2} \E H_{k_2}^2 = O \Big
(\eps^{-1}(1-\eps)^{k_1} \Big ) \, .$$ We put this into
(\ref{mom2step1}) and use (\ref{surviveprobneg}) to estimate $$\E
\Big [ (\sum _{k=r/2}^{r} H_k)^2 \mid H_{r/2}>0 \Big ] \leq O\Big
( {\eps^{-1} \over (1-\eps)^{r/2} }\Big ) \sum _{k_1=r/2}^r
\sum_{k_2=r/2}^{k_1} \eps^{-1} (1-\eps)^{k_1} = O( \eps^{-3} r )
\, ,$$
concluding our proof. \qed \\

\section{{\bf \large A lower bound on component size}}
\label{lower}

\subsection{A coupling}
Let $G$ be a transitive graph on $n$ vertices with vertex degree
$d$ and let $v \in G$ be an arbitrary vertex of $G$. We denote by
$V(G)$ and $E(G)$ the vertex and edge set of $G$, respectively.
The {\em covering tree} of $G$ rooted at $v$ is a pair $(T, \La)$
where $T$ is an infinite $d$-regular tree rooted at a vertex
$\rho$ and $\La$ is a function $\La : T \to V(G)$ satisfying:
\begin{enumerate}
\item $ \La(\rho) = v \, , \hbox{ and } \La(\rho_i) = v_i  \hbox{
for } i \in \{1,\ldots ,d\}$, where $v_1, \ldots , v_d$ are the
neighbors of $v$ in $G$ and $\rho_1, \ldots , \rho_d$ are the
children of $\rho$ in $T$.

\item For $w \in T\setminus \{\rho\}$ we have
$$ \Big \{ \La(w_1), \ldots ,
\La(w_{d-1}) \Big \} = \Big \{ v \, : \, (\La(w), v) \in E(G)
\hbox { and } v \neq \La(w^-) \Big \} \, ,$$ where $w_1, \ldots,
w_{d-1}$ are the $d-1$ children of $w$ and $w^-$ is the immediate
ancestor of $w$.
\end{enumerate}
We regard $\La$ as a labelling of the vertices of the infinite
tree $T$ with labels taking values in $V(G)$. It is clear that, up
to the choice of the arbitrary mapping between the two sets in
requirement $(2)$ of the definition, the covering tree of $G$
rooted at $v$ is unique. For the following definitions we assume
$(T, \La)$ is a covering tree of $G$ rooted at $v$ and $T_p$ is
the subgraph of $T$ obtained by performing $p$-bond-percolation on
$T$ with $p\in [0,1]$. \\

\noindent {\bf Definition 1.} A vertex $w \in T$ is called {\bf
\em impure} if there exists a vertex $u \in T$ satisfying
$$ (1)\,\,\, |u| \leq |w| \, ,\qquad \quad (2) \,\,\, \La(u) = \La(w) \, , \qquad \quad (3) \,\,\, u \lr u \wedge w \, .$$
%\begin{enumerate}
%\item $ |u| \leq |w|$, \item $\La(u) = \La(w)$, \item $\rho \lr u$.
%\end{enumerate}
See Figure $1$ for an example. A vertex $w\in T$ is called {\bf \em pure} if it is not impure. We
call a vertex $w \in T$ {\bf \em path-pure} if every vertex on
the unique path between $w$ and $\rho$ is pure. \\

\noindent {\bf Definition 2.} For a vertex set $A \subset G$, we
say $w \in T$ is {\bf \em $A$-free} if every vertex $u$ on the
unique path between $w$ and $\rho$ has $\La(u)\not \in A$. \\

%Recall that $\C(v)$ denotes the connected component containing $v$
%of $G_p$. We write $d_p(v,u)$ for the graph distance between $v$
%and $u$ in $G_p$ and we denote
%$$ B_p(v,r) = \{ u \in \C(v) : d_p(v,u) \leq r \} \, ,$$
%$$ B_p(v,r) = \{ u \in \C(v) : d_p(v,u) = r \} \, .$$
For a subset of vertices $A \subset G$ and two vertices $u,v \not
\in A$ we write $d^A_p(u,v)$ for the graph distance between $u$
and $v$ in $G_p \setminus A$ (i.e., in the graph obtained from
$G_p$ by removing the vertices of $A$ and all the edges adjacent
to $A$). We denote
$$ B^A_p(v,r) = \{ u  : d^A_p(v,u) \leq r \} \, ,$$
$$ \partial B^A_p(v,r) = \{ u  : d^A_p(v,u) = r \} \, .$$

%\noindent {\bf Definition 3.} For a vertex set $A \subset G$, we
%write $\{ x \stackrel{A^c} \lr v \}$ if $x$ is connected to $v$ in
%$G_p$ by a path which has no vertex in $A$. We write $\C(v)$ for
%the connected component containing $v$ in $G_p$ and denote by
%$\C_{A^c}(v)$ the component of $v$ in $(G\setminus A)_p$. In other
%words
%$$ \C_{A^c}(v) = \Big \{ x \, : \, x \stackrel{A^c} \lr v \Big \} \, .$$
%\\

The following random variable plays a key role in the proofs. For
an integer $r \geq 0$ and a vertex subset $A \subset G$ we define
$X^{A}_r$ by
$$ X^A_r =  \Big | \Big \{ w \in T \, : \, |w| = r \, , \, w \lr \rho \, , \,  w \hbox{ \rm is path-pure and } w
\hbox { \rm is } A\hbox{\rm -free} \Big \} \Big | \, .$$
\vspace{.1 cm}
\begin{prop}[Coupling] \label{couple} Let $G$ be a $d$-regular graph and $(T, \La)$ its covering tree. Recall the definition of
$H_r$ from Section \ref{trees}. For any $p\in [0,1]$ and any
subset $A \subset G$ there exists a coupling of $G_p$ and $T_p$
such that \be \label{coupeq} X^A_r \leq |
\partial B^A_p(v,r) | \leq H_r \, . \ee
%\Big | \Big \{ w \in T : |w| \leq k \, , \, w \lr \rho \, , \,  w
%\hbox{ \rm is path-pure and } w \hbox { \rm is } A\hbox{\rm -free}
%\Big \} \Big | \, .\ee
\end{prop}
\noindent {\bf Proof.} We recall a {\em breadth first search}
process which explores $B^A_p(v,r)$ in the spirit of Martin-L\"of
\cite{M} and Karp \cite{K}. In this exploration process, vertices
of $G$ can be either {\em explored}, {\em active} or {\em
neutral}. The active vertices are ordered in a queue which
initially contains only $v$ while the rest of the vertices are
neutral. We define a height function $h:G\setminus A \to
\mathbb{Z}^+ \cup \{\infty\}$ which is updated as the exploration
process runs. Initially we have $h(v)=0$ and $h(u)=\infty$ for any
$u \neq v$. At step $t$ of the process, we take out the first
active vertex $x_t$ from the queue. If $h(x_t)=r$ we mark it
explored, and proceed with the next step of the process. If
$h(x_t) < r$, we ``explore'' the edges between itself and its
neutral neighbors in $G\setminus A$. That is, if $x_t$ has $k$
neutral neighbors in $G\setminus A$, we examine these $k$ edges
and check whether they are open or closed (where open or closed
edges correspond to retained or deleted edges in the percolation,
respectively). For each open edge, we mark the corresponding
neutral neighbor $u$ as active, put $h(u)=h(x_t)+1$ and add $u$ to
the end of the queue. We then mark $x_t$ as explored and proceed
with the next step of the process. The process ends when the queue
empties. We make two important observations. First, the BFS tree
structure of this process guarantees that we will never have in
the queue two vertices which have height difference strictly
larger than $1$. Secondly, at the end of this process we have
$$  \partial B^A_p(v,r) = \Big \{ x \, : \, h(x) =r \Big \} \, .$$

We couple this exploration process, in a natural way, with
percolation on $T$. Each edge explored in the exploration process
will correspond to precisely one edge in the tree and these two
edges are coupled such that they are open or closed together; the
rest of the edges in the tree will be open with probability $p$
and closed otherwise, independently of all other choices. We will
also have a correspondence between the {\em vertices} of
$B^A_p(v,r)$ and the vertices of $T$. The correspondence (of both
edges and vertices) is done in the following recursive manner.
Recall that $v_1, \ldots, v_d$ are the neighbors of $v$ in $G$ and
$\rho_1, \ldots, \rho_d$ are the children of $\rho$ in $T$.
Initially, the vertex $v$ corresponds to $\rho$ and the edges
$(v,v_i)$ and $(\rho, \rho_i)$ correspond to each other and are
coupled such that they are open or closed together. At step $t>1$,
let $x_t$ be the active vertex explored and let $y_1, \ldots, y_k$
be its {\em neutral} neighbors in $G\setminus A$; note that $k
\leq d-1$. Let $w_t$ be the vertex corresponding to $x_t$ in $T$
and denote by $u_1, \ldots u_k$ the children of $w_t$ in $T$ such
that $\La(w_i) = y_i$ for $1 \leq i \leq k$. We now have the edge
$(x_t, y_i)$ and $(w_t, u_i)$ correspond to each other and we
couple such that they are open or closed together, for $1 \leq i
\leq k$. We also have the graph vertex $y_i$ correspond to the
tree vertex $u_i$ for each {\em open} edge $(x_t,y_i)$. This
finishes the description of our coupling.

Under this coupling, the upper bound of (\ref{coupeq}) is obvious
and we need to prove the lower bound. We do this in two steps.

\noindent {\bf Claim 1.} If a graph vertex $y\in G$ corresponds in
our coupling to a tree vertex $u\in T$, then $u \lr \rho$ and
$h(y) = |u|$ and $v=\La(u)$.

\noindent {\bf Proof.} Follows easily from the definition of our
coupling by induction on $h(y)$.

\noindent {\bf Claim 2.} If a tree vertex $w \in T$ is path pure
and $A$-free and has $|w|=\ell$ and $w \lr \rho$, then $\La(w)$
corresponds in our coupling to $w$ and for every pure child $w^+$
of $w$, if the tree edge $(w, w^+)$ is open, then it corresponds
in our coupling to $(\La(w), \La(w^+))$.

\noindent {\bf Proof.} We prove by induction on $\ell$. The
assertion is obvious for $\ell=0$. Let $\ell >0$ and assume $w$
satisfies the assumptions of the claim. Denote the path from
$\rho$ to $w$ in $T$ by $\rho=w_0, w_1, \ldots, w_\ell = w$. We
apply our induction hypothesis on $w_{\ell -1}$ and deduce that
$\La(w_{\ell-1})$ corresponds to $w_{\ell-1}$ and hence
$h(\La(w_{\ell-1}))=\ell-1$ by Claim $1$. We also deduce that the
tree edge $(w_{\ell-1}, w_\ell)$ corresponds to the edge
$(\La(w_{\ell-1}), \La(w_\ell))$ whence the edge
$(\La(w_{\ell-1}), \La(w_\ell))$ is open in $G_p$. Consider the
time $t$ when $\La(w_{\ell-1})$ was the active vertex $x_t$ taken
out of the queue. We claim that at that time the graph vertex
$\La(w_\ell)$ was neutral. Assume otherwise, then $\La(w_\ell)$ is
active or explored at time $t$ and since
$h(\La(w_{\ell-1}))=\ell-1$ we have that $h(\La(w_{\ell})) \leq
\ell$ by our first observation from before. We deduce that
$\La(w_\ell)$ corresponds to some tree vertex $u$ where $u \neq
w_\ell$ (because this correspondence exists at time $t$). By claim
$1$ we have that $u \lr \rho$ and $\La(u) = \La(w_\ell)$ and $|u|
\leq |w_\ell|$, hence $w_\ell$ is not pure and we have reached a
contradiction.

We learn that $\La(w_\ell)$ was neutral at time $t$ and since the
edge $(\La(w_{\ell-1}), \La(w_\ell))$ is open. Thus, it
will be examined at time $t$, concluding our proof. \qed \\

\noindent{\bf Remark.} There can be strict inequality in the lower bound of
Proposition \ref{couple}. Indeed, it may happen that $w$ is impure because of some vertex
$u$, but $u$ is impure itself (or the path between $u$ to $u\wedge w$ is not a pure path),
and in this case we may see $\La(w)$ in the exploration process on $G$, but not count it in
$X^A_r$. For example, in Figure 1, if the edge between $5$ and $4$ in the left side of the
tree was open, then both of these $4$'s would be impure.

%We conclude the proof by observing that since the tree edge
%$(w_{\ell-1},w_\ell)$ is open then so is
%$(\La(w_{\ell-1}), \La(w_\ell))$. \qed \\

\subsection{A key lemma} The following key lemma utilizes the non-backtracking
random walk to provide a lower bound on the expected size of
$B^A_p(v,r)$. In all the lemmas of the rest of Section \ref{lower}
we have a transitive graph $G$ with vertex degree $d$. Recall from
Section \ref{rw} that $\p^t(v,v)$ is the return probability after
$t$ steps of the non-backtracking random walk.

\begin{lemma} \label{keylemma} Let $v$ be a uniform random vertex of $G$. Then for any $p\in [0,1]$ and any $A \subset G$ we have
$$ \E | \partial B^A_p(v,r) | \geq (p(d-1))^r \Big [ 1 - {r|A| \over n} - {d \over
d-1} \sum _{h=2}^r \sum_{k=2}^h \sum
_{j=1}^{k-1}(p(d-1))^{k-j}\p^{h+k-2j-1}(v,v) \Big ]\, .$$
\end{lemma}
\noindent {\bf Proof.} Let $(T, \La)$ be the covering tree of $G$
rooted at $v$ and consider $T_p$. By proposition \ref{couple} it
suffices to prove the estimate of the lemma on $\E X^A_r$. To that
aim, let $W_r$ be a uniform random vertex from the set $\{ w \in T
\, : \, |w| = r \}$ and let $(W_0, W_1, \ldots, W_r)$ denote the
random path from $W_0=\rho$ to $W_r$. For triples of integers
$(j,k,h)$ satisfying $1 \leq j < k \leq h \leq r$ denote by
$X_{j,k}^{(h)}$ the random variable
$$ X_{j,k}^{(h)} = \Big | \Big \{ u \in T \, : \, |u|=k \, , \, u
\wedge W_h = W_j \, , \, \La(u) = \La(W_h) \, , \,u \lr W_j \Big
\} \Big | \, .$$ This random variable counts the number of
vertices of height $k$ in $T$ which connect to the path $(W_0,
\ldots, W_r)$ at $W_j$ and make the vertex $W_h$ impure. Indeed,
observe that by definition, if $X_{j,k}^{(h)} =0$ for all triples
with $1 \leq j < k \leq h \leq r$ then $W_r$ is path-pure. We have
that
$$ \E X^A_r = d(d-1)^{r-1} \prob \Big ( W_r \lr
\rho \, , \, W_r \hbox{ \rm is path-pure and } W_r \hbox{ \rm is }
A\hbox{\rm -free} \ \Big ) \, .$$ For any instance of $W_r$, the
random variable $X_{j,k}^{(h)}$ is determined by percolation on
edges which are {\em not} the edges of the path $W_0, \ldots ,
W_r$. Hence $X_{j,k}^{(h)}$ is independent of the event $\{W_r \lr
\rho\}$ for all triples $(j,k,h)$. Similarly, the event $\{W_r
\hbox{ \rm is } A\hbox{\rm -free}\}$ is independent of $\{W_r \lr
\rho\}$ and is implied by the event $S=0$ where $S = |\{ h \leq r
: W_h \in A\}|$.
 Hence,
$$ \E X^A_r \geq (p(d-1))^r \prob \Big ( X_{j,k}^{(h)} = 0 \hbox{
\rm for all } (j,k,h) \, , \, S = 0 \Big ) \, .$$ Since our
initial vertex $v$ was chosen uniformly at random we have $\E S =
{r|A|\over n}$. By Markov's inequality we get that \be
\label{keystep4}   \E X^A_r \geq (p(d-1))^r \Big [ 1 - {r|A| \over
n} - \sum _{h=2}^r \sum_{k=2}^h \sum _{j=1}^{k-1} \E X_{j,k}^{(h)}
\Big ] \, .\ee We are left to estimate from above $\E
X_{j,k}^{(h)}$. Let $U_{j,k}$ be an independent random uniform
vertex from the set
$$\Big \{ u \in T \, : \, |u|=~k \, , u \hbox{ is a descendant of } W_j
\Big \} \, ,$$ and note that
$$ \E X_{j,k}^{(h)} = (d-1)^{k-j} \prob \Big ( U_{j,k} \wedge W_h =
W_j\, , \, \La(U_{j,k})=\La(W_h) \, , \, U_{j,k} \lr W_j \Big ) \,
.$$ The event $\{U_{j,k} \lr W_j\}$ is independent of the event
$\{U_{j,k} \wedge W_h = W_j\, , \, \La(U_{j,k})~=~\La(W_h)~\}$,
hence \be\label{keystep1} \E X_{j,k}^{(h)} = (p(d-1))^{k-j} \prob
\Big ( U_{j,k} \wedge W_h = W_j\, , \, \La(U_{j,k}) = \La(W_h)
\Big ) \, .\ee Let $U_{j,k}^1$ be the child of $W_j$ such that
$U_{j,k}$ is a descendant of $U_{j,k}^1$ (or if $k=j+1$ we take
$U_{j,k}^1 = U_{j,k}$).  For $i,\ell \in \{1,\ldots, d-1\}$ denote
by $\A(i,\ell)$ the event
$$ \A(i,\ell) = \Big \{ W_j = w \, , W_{j+1} = w_i \, ,U_{j,k}^1 = w_\ell
\Big \} \, ,$$ where $w\in T$ has $|w|=j$ and $w_1, \ldots,
w_{d-1} \in T$ are the children of $w$ in $T$. Note that $U_{j,k}
\wedge W_h = W_{j}$ if and only if $W_{j+1} \neq U_{j,k}^1$. Hence
\be \label{keystep2} \prob \Big ( U_{j,k} \wedge W_h =
W_j,\La(U_{j,k}) = \La(W_h) \mid W_j \Big ) = { \mathop
{\displaystyle \sum _{i \neq \ell}^{d-1}} \prob \Big (
\La(U_{j,k}) = \La(W_h) \mid \A(i,\ell) \Big ) \over (d-1)^2} \, .
\ee Given that $U_{j,k}^1 = w_\ell$, we have that $\La(U_{j,k})$
is distributed as the end vertex of a non-backtracking random walk
of length $k-j-1$ on $G$ starting with the edge
$(\La(w),\La(w_\ell))$. Similarly, given that $W_{j+1} = w_i$, the
vertex $\La(W_{h})$ is distributed as the end vertex of an
independent non-backtracking random walk of length $h-j-1$ on $G$
starting with the edge $(\La(w),\La(w_i))$. We deduce that
$$ \prob \Big ( \La(U_{j,k}) = \La(W_h) \mid \A(i,\ell) \Big ) =
\sum_{y \in V} \prob _{(\La(w),\La(w_\ell))}(X_{k-j-1}^{(2)}=y)
\prob _{(\La(w),\La(w_i))}(X_{h-j-1}^{(2)}=y) \, ,$$ where
$\{X_t\}$ is the non-backtracking random walk (see the notation of
Section \ref{rw}). This together with Lemma \ref{loop} implies
that $$ \mathop {\displaystyle \sum _{i \neq \ell}^{d-1}} \prob
\Big ( \La(U_{j,k}) = \La(W_h) \mid \A(i,\ell) \Big ) \leq d(d-1)
\p^{h+k-2j-1}(v,v) \, .$$ This together with (\ref{keystep1}) and
(\ref{keystep2}) gives that $$ \E X_{j,k}^{(h)} \leq {d \over d-1}
(p(d-1))^{k-j}\p^{h+k-2j-1}(v,v) \, .$$ Putting this into
(\ref{keystep4}) completes the proof of
the lemma. \qed \\

\subsection{A second moment argument} The following two lemmas bound from below the
probability that $B^{A}_p(v, r)$ is large. The two lemmas handle
the cases $\eps>0$ and $\eps<0$, since our estimates from Section
\ref{trees} change when $\eps$ changes sign. The proofs of the two
are the same, except that we apply Lemma \ref{treecalc} in Lemma
\ref{central} and Lemma \ref{treecalcneg} in Lemma
\ref{centralneg}. To ease the reading of these lemmas, the reader
is advised to think of the three significant parameters $\eps, r$
and $M$ as taking the values $\eps \approx n^{-1/3}$, $r \approx
n^{1/3}$ and $M \approx n^{2/3}$ (scaling window width, diameter
and volume, respectively). These will be the values we will use
for the proof of Theorem \ref{mainthm1}.

\begin{lemma} \label{central} Let $\eps \in (0,1/2)$ and put $p={1 + \eps \over d-1}$.
Denote by $v$ a uniform random vertex of $G$ and let $M$ and $r$
be two integers satisfying
\begin{enumerate}
\item $\quad \displaystyle {d \over d-1} \sum _{h=2}^r
\sum_{k=2}^h \sum _{j=1}^{k-1}(1+\eps)^{k-j}\p^{h+k-2j-1}(v,v)
\leq {1 \over 2} \, ,$

%\item $\quad \displaystyle |A| \leq \ell \, ,$

%\item $\quad \displaystyle \ell r \leq n/4 \, ,$

\item $\quad \displaystyle 96 M < \eps^{-2} [(1+\eps)^r - (1+\eps)^{r/2}]
(1-e^{- \eps r/4}) \, .$
\end{enumerate}
Then there exists some fixed $c>0$ such that if $A \subset G$ has
$|A| \leq {n \over 4r}$, then we have
$$ \prob \Big ( | B^{A}_p(v, r) | \geq M  \Big
) \geq c \eps \Big (1-(1+\eps)^{-r/2}\Big )^2 \Big (1-e^{ {-r \eps
/ 4}}\Big )^3 \, .$$
\end{lemma}
\noindent {\bf Proof.} For notational convenience we write
$\partial B_k$ for $|\partial B^{A}_p(v,k)|$ so that $|B^{A}_p(v,
r)|=\sum_{k=0}^r \partial B_k$. We have
$$ \prob \Big ( \Big | B^{A}_p(v, r) \Big | \geq M \Big
) \geq  \prob \Big ( \partial B_{r/2} > 0 \Big ) \prob \Big (
\sum_{k=r/2}^r \partial B_k \geq M \mid \partial B_{r/2} > 0 \Big)
\, .$$ By Cauchy-Schwartz we get \be\label{comp1} \prob \Big (
\Big | B^{A}_p(v, r) \Big | \geq M \Big ) \geq { \Big [ \E
\partial B_{r/2} \Big ]^2 \over \E \Big [ (\partial B_{r/2}) ^2 \Big ]} \, \prob \Big ( \sum_{k={r/2}}^r
\partial B_k \geq M \mid \partial B_{r/2} > 0 \Big) \, .\ee
By definition $\E[\partial B_k \mid \partial B_{r/2} > 0] = {\E
\partial B_k \over \prob (\partial B_{r/2}~>~0~)}$ for any
$k\geq r/2$. Together with Proposition \ref{couple} this implies
that \be\label{comp2} \E \Big [ \sum_{k= {r/2}}^r \partial B_k
\mid
\partial B_{r/2} > 0\Big ] \geq { \sum_{k={r/2}}^r \E \partial B_k \over \prob (H_{{r/2}} >0 )} \, .\ee
Lemma \ref{keylemma} together with assumption $(1)$ and our
assumption on $A$ imply that for any $k \leq r$ we have
\be\label{comp4} \E
\partial B_k \geq (1+\eps)^k \Big [ {1 \over 2} - {r |A| \over n}
\Big ] \geq { (1+\eps)^k \over 4 } \, .\ee We put this in
(\ref{comp2}) which together with (\ref{surviveprob}) of Lemma
\ref{treecalc} yields $$ \E \Big [ \sum_{k= {r/2}}^r \partial B_k
\mid \partial B_{r/2} > 0 \Big ] \geq {1 \over 48} \eps^{-2}
[(1+\eps)^r - (1+\eps)^{r/2}] (1-e^{- \eps r/4}) \, .$$ Assumption
$(2)$ and the previous display imply that
$$M < {1 \over 2} \E \Big [ \sum_{k= {r/2}}^r \partial B_k \mid
\partial B_{r/2} > 0 \Big ] \, .$$
Hence we can use the estimate $\prob (Z > y) \geq (\E Z - y)^2/\E
Z^2$, valid for any non-negative random variable $Z$ and $y < \E
Z$ (this estimate follows easily from Cauchy-Schwartz), and get
\be \label{comp3} \prob \Big ( \sum_{k={r/2}}^r
\partial B_k \geq M \mid \partial B_{r/2}
> 0 \Big) \geq { c \eps^{-4} [(1+\eps)^r - (1+\eps)^{r/2}]^2
(1-e^{- \eps r/4})^2 \over \E \Big [ (\sum_{k={r/2}}^r
\partial B_k)^2 \mid \partial B_{r/2} > 0\Big ] } \, ,\ee
for some $c>0$. Put $\A = \{\partial B_{r/2} > 0\}$ and $\B =
\{H_{r/2}>0\}$ and note that Proposition \ref{couple} allows us to
couple such that $\A \subset \B$. Since $\A \subset \B$, for any
non-negative random variable $X$ we have
$$\E[X \mid \A] = {\E [ X {\bf 1}_{\A} ] \over \prob(\A)} \leq
{\prob (\B) \over \prob (\A)} \, \E[ X \mid B]  \, .$$ We put
$X=(\sum_{k={0}}^r \partial B_k)^2$ and use Proposition
\ref{couple}, which allows us to couple such that $\partial B_k
\leq H_k$, to get \be \label{uselater2} \E \Big [
(\sum_{k={r/2}}^r
\partial B_k)^2 \mid \partial B_{r/2} > 0\Big ] \leq { \prob (H_{r/2}>0) \over\prob ( \partial B_{r/2}
> 0 )} \, \E \Big [ (\sum_{k= {r/2}}^r
H_k)^2 \mid H_{r/2} > 0 \Big ] \, .\ee We use Cauchy-Schwartz to
bound from below the denominator and Lemma \ref{treecalc} to bound from above the two parts of the numerator
$$ \E \Big [ (\sum_{k={r/2}}^r
\partial B_k)^2 \mid \partial B_{r/2} > 0\Big ] \leq { \E \Big [ (\partial B_{r/2} )^2 \Big ] \over \Big [ \E \partial B_{r/2} \Big
]^2 } O \Big ( \eps^{-3} (1+\eps)^{2r} (1-e^{\eps r/4})^{-1}\Big )
\, .$$ We put this in (\ref{comp3}) and get that
$$ \prob \Big
( \sum_{k={r/2}}^r \partial B_k \geq M \mid \partial B_{r/2} > 0
\Big) \geq c { \Big [ \E \partial B_{r/2} \Big ]^2 \over \E \Big [
(\partial B_{r/2}) ^2 \Big ] } \eps^{-1} [1 - (1+\eps)^{-r/2}]^2
(1-e^{- \eps r/4})^3 \, ,$$ for some $c>0$. We plug this into
(\ref{comp1}) and get \be\label{comp5} \prob \Big ( \Big |
B^{A}_p(v, r) \Big | \geq M \Big ) \geq c { \Big [ \E
\partial B_{r/2} \Big ]^4 \over \E \Big [ (\partial B_{r/2}) ^2 \Big ]^2  } \eps^{-1}
[1 - (1+\eps)^{-r/2}]^2 (1-e^{- \eps r/4})^3 \, .\ee By
(\ref{treemom2}) of Lemma \ref{treecalc} and Proposition
\ref{couple} we have
$$ \E \Big [ (\partial B_{r/2}) ^2 \Big ] \leq O \Big ( \eps^{-1} (1+\eps)^{r} \Big ) \, , $$
and putting $k=r/2$ in (\ref{comp4}) and the result into
(\ref{comp5}) yields the assertion of the lemma. \qed \\

\begin{lemma} \label{centralneg} Let $\eps \in (0,1/2)$ and put $p={1 - \eps \over d-1}$.
Denote by $v$ a uniform random vertex of $G$ and let $M$ and $r$
be two integers satisfying
\begin{enumerate}
\item $\quad \displaystyle {d \over d-1} \sum _{h=2}^r
\sum_{k=2}^h \sum _{j=1}^{k-1}(1-\eps)^{k-j}\p^{h+k-2j-1}(v,v)
\leq {1 \over 2} \, ,$

\item $\quad \displaystyle 96 M <  \eps^{-1} r [ (1-\eps)^{r/2} -
(1-\eps)^{r} ] \, .$
\end{enumerate}
Then there exists some fixed $c>0$ such that if $A \subset G$ has
$r|A| \leq n/4$ then we have
$$ \prob \Big ( | B^{A}_p(v, r) | \geq M  \Big
) \geq c \eps^3 r^2 (1-\eps)^{r} [ (1-\eps)^{r/2} - (1-\eps)^{r}
]^2 \, .$$
\end{lemma}

\noindent {\bf Proof.} The proof carries on precisely as in the
previous lemma up to (\ref{comp4}). Instead we use assumption
$(1)$ and Lemma \ref{keylemma} to estimate \be\label{comp4neg} \E
\partial B_k \geq { (1-\eps)^k \over 4 } \, ,\ee which together
with (\ref{surviveprobneg}) of Lemma \ref{treecalcneg} yields
$$ \E \Big [ \sum_{k= {r/2}}^r \partial B_k
\mid \partial B_{r/2} > 0 \Big ] \geq {1 \over 48} \eps^{-1} r [
(1-\eps)^{r/2} - (1-\eps)^{r} ]  \, .$$ As before, by assumption
$(2)$ we deduce that for some $c>0$
 \be \label{comp3neg} \prob \Big ( \sum_{k={r/2}}^r
\partial B_k \geq M \mid \partial B_{r/2}
> 0 \Big) \geq { c \eps^{-2} r^2 [(1-\eps)^{r/2} -
(1-\eps)^{r} ]^2 \over \E \Big [ (\sum_{k={r/2}}^r
\partial B_k)^2 \mid \partial B_{r/2} > 0\Big ] } \, .\ee
Inequality (\ref{uselater2}) still holds and we use it to estimate
the denominator of the last display. As before, we use
Cauchy-Schwartz to bound the denominator of (\ref{uselater2}) and
Lemma \ref{treecalcneg} to bound the two parts of the numerator of
(\ref{uselater2}). This gives
$$ \E \Big [ (\sum_{k={r/2}}^r
\partial B_k)^2 \mid \partial B_{r/2} > 0\Big ] \leq { \E \Big [ (\partial B_{r/2} )^2 \Big ] \over \Big [ \E \partial B_{r/2} \Big
]^2 } O \Big ( \eps^{-3}\Big ) \, .$$ We put this in
(\ref{comp3neg}) and the result into (\ref{comp1}) to get that for
some $c>0$ we have \be\label{comp5neg} \prob \Big ( \Big |
B^{A}_p(v, r) \Big | \geq M \Big ) \geq c \eps r^2  { \Big [ \E
\partial B_{r/2} \Big ]^4 \over \E \Big [ (\partial B_{r/2}) ^2 \Big ]^2  } [ (1-\eps)^{r/2} -
(1-\eps)^{r} ]^2 \, .\ee By Proposition \ref{couple} and
(\ref{treemom2neg}) of Lemma \ref{treecalcneg} we get
$$ \E \Big [ (\partial B_{r/2}) ^2 \Big ] \leq O \Big (\eps^{-1}
(1-\eps)^{r/2} \Big ) \, ,$$ and putting $k=r/2$ in
(\ref{comp4neg}) and that into (\ref{comp5neg}) gives that
$$ \prob \Big ( \Big | B^{A}_p(v, r) \Big | \geq M \Big ) \geq c
\eps^3 r (1-\eps)^{r} [ (1-\eps)^{r/2} - (1-\eps)^{r} ]^2 \, .$$
\qed

\section{{\bf \large Proof of Theorems \ref{mainthm1} and \ref{mainthm2}}}

To prove Theorems \ref{mainthm1} and \ref{mainthm2}, we employ a
process which explores neighborhoods of vertices in $G_p$. For a
fixed number $r$, to be chosen later, the process explores
neighborhoods of randomly chosen vertices up to distance $r$,
excluding the vertices it has seen before. For the following,
recall the definition of $B_p^A(v, r)$ from Section \ref{lower}.
The process runs as follows. We start by choosing a uniform random
vertex $v_1$ and putting $V_1 = B_p(v_1, r)$. At each step $t \geq
2$ we choose a random uniform vertex $v_t$ of $G$ and put $V_t =
V_{t-1} \cup B^{V_{t-1}}_p(v_t, r)$. The process ends when we
exhaust all the vertices in the graph, i.e., when $|V_t| = n$. For
convenience we write $V_0 = \emptyset$.

For some fixed number $M$ we wish to study if this process has
encountered a neighborhood of size at least $M$. We introduce some
notation. Write $I_t$ for the indicator random variable for the
event $\{ |B^{V_{t-1}}_p (v_t, r)| > M \}$. It is clear that if
there exists $t$ with $I_t=1$ then $|\C_1|\geq M$. Hence \be
\label{final} \prob \Big ( |\C_1| \geq M \Big ) \geq \prob \Big
(\exists t \hbox{ {\rm with} } I_t = 1 \Big ) \, .\ee

% and write $S_t = \sum _{j=1}^t I_j$. Write $\tau$ for the
%stopping time
%$$ \tau = \min \{ t : V_t \geq M \} \, .$$
%Observe that if $I_t = 1$ then we must have $\tau \leq t$; we
%deduce that $S_\tau \in \{0,1\}$. It is clear that if there exists
%$t$ with $I_t=1$ then $|\C_1|\geq M$ and in particular $S_\tau =
%1$ implies $|\C_1|\geq M$. Thus, \be \label{final} \prob \Big (
%|\C_1| \geq M \Big ) \geq \prob (S_\tau \geq 1) = \E S_\tau \,
%.\ee

\begin{lemma} \label{volbound} Let $p={1 + \eps \over d-1}$ and
$T>0$. If $\eps >0$, then we have
$$ \E |V_T| \leq {2 T (1+\eps)^{r+1} \over \eps} \, ,$$
and if $\eps < 0$, then we have
$$ \E |V_T| \leq 2 T r \, .$$
\end{lemma}
\noindent {\bf Proof.} We have that
$$ |V_T| = \sum_{t=1}^T |B^{V_{t-1}}_p(v_t,r)| \, .$$
Once we condition on $V_{t-1}$ and $v_t$, Proposition \ref{couple}
allows us to couple such that $|\partial B^{V_{t-1}}_p(v_t,k)|
\leq~H_k$. Hence, for any $t \leq T$ and $\eps>0$ we have
$$ \E |B^{V_{t-1}}_p(v_t,r)| \leq \sum _{k=0}^r \E H_k \leq
{d \over d-1} {(1+\eps)^{r+1} \over \eps } \, ,$$ which concludes the proof for
$\eps>0$ (since $d\geq 3$). The proof for $\eps<0$ goes similarly by bounding $\E
|B^{V_{t-1}}_p(v_t,r)| \leq 2r$.
\qed \\

\noindent {\bf Proof of Theorem \ref{mainthm1}.} Recall that $p =
{1 + \lambda n^{-1/3} \over d-1}$ for some fixed $\lambda \in
\mathbb{R}$. Since $\prob ( |\C_1| \geq M )$ is increasing with
$p$ and $\prob ( |\C_1| \leq M )$ is decreasing with $p$ we may
assume that $|\lambda| \geq 1$ and the result follows for $|\lambda|<1$. Theorem
$1.2$ of \cite{NP3} (or Proposition $1$ of \cite{NP2}) states that
for any $\alpha>0$ there exists $A=A(\alpha, \lambda)>0$ such that
for {\em any} graph with maximum degree $d \in [3,n-1]$ we have
$$ \prob \Big ( |\C_1| > An^{2/3} \Big ) \leq \alpha \, .$$
Thus the upper bound on $|\C_1|$ implied in Theorem \ref{mainthm1}
is already proved without the need of condition (\ref{maincond1}).
\\

For the lower bound, fix some small $\delta = \delta(\lambda)>0$
and $\gamma=\gamma(\lambda) \in (0,1)$, to be chosen later. We put
$\eps(n) = \lambda n^{-1/3}$ and $r = \gamma  n^{1/3}$ and $M =
\delta n^{2/3}$. Let $h \leq r$ and $t \leq 2h$ be a two positive
integers. If $(j,k)$ is a pair of integers satisfying $2\leq k \leq h$ and
$1 \leq j \leq k-1$ and $h+k-2j-1=t$ then we must have $h-t+1 \leq
k \leq h$. We deduce that the number of pairs $(j,k)$ satisfying
the above requirements is at most $t$. Thus, we bound
\begin{eqnarray*} \sum _{h=2}^r \sum_{k=2}^h \sum
_{j=1}^{k-1}(p(d-1))^{k-j}\p^{h+k-2j-1}(v,v) &\leq& e^{|\lambda|\gamma}
\sum _{h=2}^r \sum _{t=1}^{2h} t \p^t(x,x) \\ &\leq& e^{|\lambda|
\gamma} \gamma n^{1/3} \sum_{t=1}^{2\gamma n^{1/3}} t \p^t(x,x) \,
, \end{eqnarray*} where in the first inequality we bounded $(p(d-1))^{k-j} \leq (1+|\lambda|n^{-1/3})^r$.
Condition (\ref{maincond1}) and the last display
imply that we can choose $\gamma\in (0,1)$ small enough such that
$$ {d \over d-1} \sum _{h=2}^r \sum_{k=2}^h \sum
_{j=1}^{k-1}(1+\eps)^{k-j}\p^{h+k-2j-1}(v,v) \leq {1 \over 2} \,
.$$ For any $\gamma>0$ we can choose $\delta>0$ so small such that
the two assumptions of Lemma \ref{central} or Lemma
\ref{centralneg} (according to whether $\lambda >0$ or $\lambda <
0$) hold. We deduce that there exists a constant $c =
c(\lambda)>0$ such that \be \label{point} \prob \Big ( |
B^{V_{t-1}}_p(v_{t}, r) | \geq M \, \Big | \, |V_{t-1}| \leq {n
\over 4r} \Big ) \geq c n^{-1/3} \, .\ee For some positive integer
$T$ denote by $\A_T$ the event
$$ \A_T = \Big \{ \forall t \leq T \, \, I_t = 0  \hbox{ \rm{and} }
V_T \leq {n \over 4r} \Big \} \, .$$ We prove by induction that
$\prob (\A_T) \leq (1-cn^{-1/3})^T$. Indeed for $T=1$ it is
obvious by (\ref{point}). We have
$$ \prob ( \A_T ) = \prob \Big ( \A_{T-1} \hbox{ {\rm and} } | B^{V_{T-1}}_p(v_{T}, r) | < M  \hbox{ {\rm and}
} | B^{V_{T-1}}_p(v_{T}, r) | < {n \over 4r} - |V_{T-1}| \Big ) \,
.$$ Hence we can bound
$$ \prob ( \A_T ) \leq \prob( \A_{T-1} ) \prob \Big ( | B^{V_{T-1}}_p(v_{T}, r)
| < M \, \Big | \, \A_{T-1} \Big ) \leq \prob( \A_{T-1} )
(1-cn^{-1/3}) \, ,$$ where the last inequality is done by
conditioning on the sets $V_1, V_2, \ldots, V_{T-1}$ and using
(\ref{point}). We now have
$$ \prob \Big (  \forall t \leq T \, \, I_t = 0 \Big ) \leq \prob (
\A_T ) + \prob \Big ( |V_T| \geq {n \over 4r} \Big ) \leq (1-
c n^{-1/3})^T + { 4 r \E V_T \over n } \, ,$$ where the
last inequality follows from our estimate on $\prob (\A_T)$ and
Markov's inequality. Put $T= \alpha n^{1/3}$ for some $\alpha>0$
and use Lemma \ref{volbound} to get
$$ \prob \Big (  \forall t \leq T \, \, I_t = 0 \Big ) \leq
e^{-c \alpha} +  { 8\alpha \gamma e^{\lambda \gamma} \over
\lambda } \, ,$$ for $\lambda > 0$ and
$$ \prob \Big (  \forall t \leq T \, \, I_t = 0 \Big ) \leq
e^{-c \alpha} + {8\gamma^2 \alpha} \, ,$$ for $\lambda < 0$.
Putting $\alpha = \gamma^{-1/2}$ yields that we can make both of
these probabilities sufficiently small by taking $\gamma$
small enough. This together with (\ref{final}) concludes the proof.  \qed \\

\noindent {\bf Proof of Theorem \ref{mainthm2}.} We proceed
similarly to the proof of Theorem \ref{mainthm1}. Take $r$ defined
in Theorem \ref{mainthm2}. Let $(h,t_1,t_2)$ be a triple of positive integers
satisfying
$$ h \leq r \, , \qquad t_1 \leq 2h \, , \qquad t_2 \leq h \wedge
t_1 \, .$$ The number of pairs $(j,k)$ satisfying $h+k-2j-1=t_1$
 and $k-j=t_2$ is at most $1$. Therefore we can bound
\begin{eqnarray*} \sum _{h=2}^r \sum_{k=2}^h \sum
_{j=1}^{k-1}(1+\eps)^{k-j}\p^{h+k-2j-1}(v,v) &\leq& \sum _{h=2}^r
\sum_{t_1=1}^{2h} \sum _{t_2=1}^{h \wedge t_1} (1+ \eps)^{t_2} \p^{t_1}(v,v) \\
&\leq& \eps^{-1} r \sum_{t=1}^{2r} [(1+\eps)^{t\wedge r}-1]
\p^{t}(v,v) \, .
\end{eqnarray*}
Thus, condition (\ref{maincond2}) implies that for large enough
$n$ we have
$$ {d \over d-1} \sum _{h=2}^r \sum_{k=2}^h \sum
_{j=1}^{k-1}(1+\eps)^{k-j}\p^{h+k-2j-1}(v,v) \leq {1 \over 2} \,
.$$ Now put $M={\delta n \eps \over \log^3(n\eps^3) }$ and observe
that if $\delta>0$ is small enough, our choices of $M$ and $r$
satisfy the two assumptions of Lemma \ref{central}. Hence, there
exists some constant $c>0$ such that \be \label{point2} \prob \Big
( | B^{V_{t-1}}_p(v_{t}, r) | \geq M \, \Big | \, |V_{t-1}| < {n
\over 4r} \Big ) \geq c \eps \, .\ee We proceed similarly to
before.  For some positive integer $T$ denote by $\A_T$ the event
$$ \A_T = \Big \{ \forall t \leq T \, \, I_t = 0  \hbox{ \rm{and} }
V_T \leq {n \over 4r} \Big \} \, .$$ We prove by induction as
before that $\prob (\A_T) \leq (1-c\eps)^T$ and deduce, using
Lemma \ref{volbound} as before, that
$$ \prob \Big (  \forall t \leq T \, \, I_t = 0 \Big ) \leq (1
-c\eps)^T + { 8 T r (1+\eps)^{r} \over \eps n } \leq e^{-c\eps T}
+ {8 T \eps \over \log^2 (n \eps^3)}  \, .$$ Taking $T= \eps^{-1}
 \log(n\eps^3)$ makes this probability tend to $0$ and
(\ref{final}) concludes the proof.  \qed \\

\section{ {\bf \large Concluding remarks and open problems}}\label{conc}

\begin{enumerate}

\item[1.] The proof of Theorem \ref{mainthm2} shows in fact that
condition (\ref{maincond2}) could be replaced by the slightly
weaker condition \be \label{maincond2alt} \limsup _{n} \,
\eps^{-1} r \sum _{t=1}^{2r} [(1+\eps)^{t\wedge r} -1] \p^t (v,v)
< \limsup_n {d(n)-1 \over d(n)} \, . \ee We have not found,
however, examples in which this condition holds and
(\ref{maincond2}) does not.

\item[2.] In the case that $G$ is the hamming $m$-cube
$\{0,1\}^m$, an upper bound on the size of the scaling
window of order $n^{-1/\log^{2/3}(n)}$ is given in \cite{BCHSS3}.
It is broadly believed (and conjectured in \cite{BCHSS3}) that the
scaling window of the hypercube is of order $\Theta(n^{-1/3})$,
however, this is still wide open.

\item[3.] The case of the high-dimension discrete torus $[m]^d$
(for large fixed $d$ and $m~\to~\infty$~) seems even harder. In
this case there is no sub-constant upper bound on the size of the
scaling-window.
% A recent work of Heydenreich and van der Hofstad
% shows the intriguing fact, that up to some logarithmic
% corrections, a scaling window of width {\em at least} $n^{-1/3}$
% occurs around $p_c(\mathbb{Z}_d)$ (i.e., the critical probability
% for bond-percolation on the infinite lattice $\mathbb{Z}_d$). \\

\item[4.] Another problem is to derive the existence of a scaling
window for expanders of low girth, where the critical probability
is not ${1 \over d-1}$. In particular, we conjecture that the
conclusions of Theorem \ref{expanderthm} hold for {\em any}
expander family, only around a different $p_c$, larger than ${1
\over d-1}$.

%\item[5.] It is easy to construct examples of transitive graphs
%$\{G_n\}$ adhering to conditions (\ref{maincond1}) and
%(\ref{maincond2}) in which the size of the largest component
%outside of the scaling window is not concentrated. For instance,
%take two disjoint complete graphs on $n$ vertices and connect them
%with a matching. It is plausible that conditions (\ref{maincond1})
%and (\ref{maincond2}) together with expansion imply the
%concentration of the size of the largest component.

\end{enumerate}

\section*{ {\bf  \large Acknowledgments}}
I am greatly indebted to Yuval Peres for his careful guidance, his
encouragement to work on this problem and for many insightful
conversations and suggestions. I warmly thank Itai Benjamini and
Gady Kozma for many inspiring conversations that eventually led to
this research. I am also grateful to Eyal Lubetzky for useful
conversations and for assisting me with non-backtracking random
walk issues, to Noga Alon for useful comments regarding expander
graphs, and to Jian Ding for correcting some mistakes in an
earlier version of this manuscript.

I also thank the Theory Group in Microsoft Research, where parts
of this research was conducted, for their kind hospitality and
support.

\medskip \noindent
{\bf Asaf Nachmias}: \texttt{asafn(at)microsoft.com} \\
Microsoft Research,
One Microsoft way,\\
Redmond, WA 98052-6399, USA.

\end{document}